\documentclass{article}
\usepackage[top=30truemm,bottom=30truemm,left=25truemm,right=25truemm]{geometry}
\usepackage{amsmath,amssymb,amsthm}
\usepackage{amscd}
\usepackage{tikz}
\usepackage{enumerate}
\usepackage{url}
\usepackage{hyperref}

\newtheorem{theorem}{Theorem}[section]

\theoremstyle{definition} %%%
\newtheorem{df}[theorem]{Definition}
\newtheorem{ex}[theorem]{Example}
\newtheorem{rem}[theorem]{Remark}

\theoremstyle{remark}
\newtheorem{pf}{Proof\!}

\theoremstyle{theorem} %%%
\newtheorem{thm}[theorem]{Theorem}
\newtheorem{lem}[theorem]{Lemma}
\newtheorem{prop}[theorem]{Proposition}
\newtheorem{cor}[theorem]{Corollary}
\newtheorem{conj}[theorem]{Conjecture}

\makeatletter
    
    \@addtoreset{equation}{section}
\makeatother

\def\C{\mathcal{C}}

\def\P{\mathcal{P}}

\def\c{\mathbf{c}}

\def\Z{\mathbb{Z}}
\def\h{\mathfrak{h}}

\def\<{\langle}
\def\>{\rangle}

\def\E{\mathcal{E}}

\def\alert{\textbf}
\def\qed{\hfill$\Box$}
\def\fin{\hfill$\Diamond$}

\def\tilde{\widetilde}

\def\P{\mathcal{P}}
\def\Z{{\mathbb{Z}}}
\def\D{\mathcal{E}}

\def\cdiag{{\theta}}
\def\ofil{\eta}

\def\C{\mathcal{C}}

\def\lm{{\lambda}}

\def\vep{{\varepsilon}}
\def\<{\langle}
\def\>{\rangle}
\def\path{\mathfrak{p}}

\def\peri{{\omega}}

\newcommand\cyl[1]{{\mathring{#1}}}
\newcommand\cy[1]{{\mathring{#1}}}
\newcommand\perid[1]{{\hat{#1}}}
\newcommand\semid[1]{{\boldsymbol{#1}}}

\newcommand\pt[1]{{p_{#1}}}

\def\LL{{\cyl{\mathcal L}}}
\def\LP{{\mathcal L}}
\def\LE{{\mathrm{RST}}}
\def\lex{{\varepsilon}}

\begin{document}
\title{On Hook Formulas for Cylindric Skew Diagrams}
\author{Takeshi Suzuki%
\footnote{Department of Mathematics, Okayama University, Japan
\newline E-mail: suzuki@math.okayama-u.ac.jp}
\ and
Yoshitaka Toyosawa%
\footnote{Graduate School of Natural Science and Technology, Okayama University, Japan
\newline E-mail: prkr5rq9@s.okayama-u.ac.jp}}
\date{}
\maketitle

\begin{abstract}
We present a conjectual hook formula concerning  the number of the standard tableaux
on "cylindric" skew diagrams.
Our formula can be seen as an extension of Naruse's hook formula for skew diagrams.
Moreover, we  prove our conjecture in some special cases.
\end{abstract}

%%%%%%%%%%%%%%%%%%%%%%%%%%%%%%%%%%%%%%%%%%%%%%
\section{Introduction}

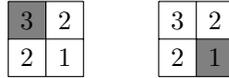
\begin{figure}[h]
\begin{center}
\begin{tikzpicture}[scale=0.5]
\fill[gray] (0,1) rectangle +(1,1);
\draw (0,0) grid +(2,2);
\foreach \a/\b/\h in {0/0/2, 0/1/3, 1/0/1, 1/1/2}
\node at (\a+0.5,\b+0.5) {$\h$};

\begin{scope}[xshift=4cm]
\fill[gray] (1,0) rectangle +(1,1);
\draw (0,0) grid +(2,2);
\foreach \a/\b/\h in {0/0/2, 0/1/3, 1/0/1, 1/1/2}
\node at (\a+0.5,\b+0.5) {$\h$};
\end{scope}
\end{tikzpicture}
\end{center}
\caption{
The excited diagrams of $\mu=(1,0)$ in $\lm=(2,2)$.
Here, the number in each cell expresses the hook length.
}\label{fig:exNaruse}
\end{figure}

The hook formula gives the number of the standard tableaux on 
Young diagrams and it was discovered in 1950's 
%by Frame, Robinson and Thrall 
\cite{FRT}.
A generalization of the hook formula 
to skew diagrams was obtained relatively recently in \cite{Nar},
where Naruse gave the following formula by 
introducing excited diagrams.
%%%%%%%%%%%%%%%%%%%%%%%%%%%%%%%%%%%%%%%%%
\begin{thm}[Naruse \cite{Nar}]%\label{thm;NHF}
Let $\lambda$ and $\mu$ be partitions with $\lambda \supset \mu$ and $|\lambda/\mu|=n$.
Then the number $f^{\lambda/\mu}$ of standard tableaux on the 
skew diagram $\lambda/\mu$ is given by
\begin{equation}\label{eq;NHF}
f^{\lambda/\mu} = n! \sum_{D \in \mathcal{E}_\lambda(\mu)} 
\prod_{x \in \lambda \setminus D} \frac{1}{h_\lambda(x)},
\end{equation}
where $\mathcal{E}_\lambda(\mu)$ denotes the set of all 
excited diagrams of $\mu$ in $\lambda$, and $h_\lm(x)$ denotes
the hook length at $x$.
\end{thm}
%%%%%%%%%%%%%%%%%%%%%%%%%%%%%%%%%%%%%%%%%
%
For example, for the partitions 
$\lm=(2,2)$ and $\mu=(1,0)$, the formula leads
$$2=3!\left(
\frac{1}{2\cdot 2\cdot 1}+\frac{1}{3\cdot 2\cdot 2}
\right)$$
(See Figure \ref{fig:exNaruse}.)
%%%%%%%%%%%%%%%%%%%%%%%%%%%%%%%%%%%%%%%%%%%%%%%%
Several proofs and generalization have been known.
Morales, Pak and Panova gave a q-analogue of the skew hook formula (\cite{MPP}).
Naruse and Okada generalized 
the skew hook formula to the case where Young diagrams are replaced by general
$d$-complete posets (\cite{NO}).

In this paper, we will treat periodic or cylindric analogue of
skew diagrams (\cite{GK,Pos}) and standard tableaux on them.

Let $\peri\in \Z_{\geqq 1}\times\Z_{\leqq -1}$.
A periodic skew diagram of period $\peri$ 
is a skew diagram consisting of infinitely 
many cells which is invariant under the parallel translation by $\omega$.
We will define a standard tableau on a periodic skew diagram
as a periodic array of natural numbers
 whose  entries increase in row and column directions.
 (See Section 3 for precise definition.)

Figure \ref{fig;intro} indicates the periodic diagram 
$$\perid{\lm}/\perid{\mu}=\lm/\mu+\Z\peri=
\{u+k\peri\mid u\in\lm/\mu,\ k\in\Z\}$$
 of period $\peri=(2,-2)$ associated with the partitions $\lm =(3,1), \mu=(0,0)$,
 and two standard tableaux on it. 
(In this case, these two tableaux exhaust all the periodic standard tableaux.)

%%%%%%%%%%%%%%%%%%%%%%%%%%%%%%%%%%%%%%%
\begin{figure}[h]%\label{fig;intro}
\begin{center}
\begin{tikzpicture}[scale=0.5]

\fill[pink] (0,0) -| (3,-1) -| (1,-2) -| (0,0) -- cycle;
\draw[ultra thick] (0,0) -| (3,-1) -| (1,-2) -| (0,0) -- cycle;
\draw (0,-1) -| (1,0);
\draw (2,-1) -- (2,0);
\foreach \x/\y/\t in {0/0/1, 1/0/2, 2/0/4, 0/-1/3}
\node at (\x+0.5,\y-0.5) {$\t$};

\begin{scope}[xshift=2cm,yshift=2cm]
\draw (0,0) -| (3,-1) -| (1,-2) -| (0,0) -- cycle;
\foreach \x/\y/\t in {0/0/1, 1/0/2, 2/0/4, 0/-1/3}
\node at (\x+0.5,\y-0.5) {$\t$};
\draw (0,-1) -| (1,0);
\draw (2,-1) -- (2,0);

\draw[thick, dotted] (2,1) -- +(3,0);
\end{scope}

\begin{scope}[xshift=-2cm,yshift=-2cm]
\draw (0,0) -| (3,-1) -| (1,-2) -| (0,0) -- cycle;
\foreach \x/\y/\t in {0/0/1, 1/0/2, 2/0/4, 0/-1/3}
\node at (\x+0.5,\y-0.5) {$\t$};
\draw (0,-1) -| (1,0);
\draw (2,-1) -- (2,0);
\draw[thick, dotted] (-2, -3) -- +(3,0);
\draw[thick, dotted] (-3, -4) -- +(3,0);
\end{scope}

\begin{scope}[xshift=10cm]
\fill[pink] (0,0) -| (3,-1) -| (1,-2) -| (0,0) -- cycle;
\draw[ultra thick] (0,0) -| (3,-1) -| (1,-2) -| (0,0) -- cycle;
\foreach \x/\y/\t in {0/0/1, 1/0/3, 2/0/4, 0/-1/2}
\node at (\x+0.5,\y-0.5) {$\t$};
\draw (0,-1) -| (1,0);
\draw (2,-1) -- (2,0);

\begin{scope}[xshift=2cm,yshift=2cm]
\draw (0,0) -| (3,-1) -| (1,-2) -| (0,0) -- cycle;
\foreach \x/\y/\t in {0/0/1, 1/0/3, 2/0/4, 0/-1/2}
\node at (\x+0.5,\y-0.5) {$\t$};
\draw[thick, dotted] (2,1) -- +(3,0);
\draw (0,-1) -| (1,0);
\draw (2,-1) -- (2,0);
\end{scope}

\begin{scope}[xshift=-2cm,yshift=-2cm]
\draw (0,0) -| (3,-1) -| (1,-2) -| (0,0) -- cycle;
\foreach \x/\y/\t in {0/0/1, 1/0/3, 2/0/4, 0/-1/2}
\node at (\x+0.5,\y-0.5) {$\t$};

\draw[thick, dotted] (-2, -3) -- +(3,0);
\draw[thick, dotted] (-3, -4) -- +(3,0);
\draw (0,-1) -| (1,0);
\draw (2,-1) -- (2,0);
\end{scope}
\end{scope}

\end{tikzpicture}
\end{center}
\caption{}\label{fig;intro}
\end{figure}
%%%%%%%%%%%%%%%%%%%%%%%%%%%%%%%%%%%%%%%%%%%%%%%%%%%%%%%%%%

The image 
of  a periodic skew diagram of period $\peri$ under 
the projection $\pi:\Z^2\to \Z^2/\Z\peri$ is called a cylindric skew diagram.
The set of standard tableaux on a periodic skew diagram can be identified 
with the set of  standard tableaux on the corresponding cylindric 
skew diagram.

We remark that the cylinder $\Z^2/\Z\peri$ has a poset structure
 induced from that of $\Z^2$, and cylindric
skew diagrams can be seen as 
$d$-complete posets consisting of infinitely many cells  (cf. \cite{Str}).

We also note  that periodic/cylindric skew diagrams parameterize a
certain class of irreducible modules over 
the Cherednik algebras (double affine Hecke algebras)
(\cite{SV,Su;ratandtrig})
and the (degenerate) affine Hecke algebras \cite{Kle, Ruff}) of type $A$,
and cylindric standard tableaux also appear in those theories.

We will introduce excited diagrams for periodic skew diagrams,
and present a conjectual hook 
formula (Conjecture \ref{conj2}) 
 concerning the number of the periodic standard tableaux
on a periodic skew diagram.
The formula in Conjecture \ref{conj2} looks similar to Naruse's 
skew hook formula,
but  for periodic skew diagrams, there are infinitely many
 excited diagrams in general and the right hand side is an infinite sum.
For example, in the case where 
$\lm=(2),\mu=(0)$ and $m=\ell=1$, our hook formula leads
$$1=2!\left(
\frac{1}{1\cdot 3}+\frac{1}{3\cdot 5}+\frac{1}{5\cdot 7}+\cdots\right)
$$
We will prove that our conjecture is correct in the following cases
(Theorem \ref{thm1} and \ref{thm2}):
\begin{itemize}
%$\bullet$ 
\item (bar case) $\lm=(n),\ \mu=(0)$ and $\peri=(1,-\ell)$.

%$\bullet$
\item (hook case) $\lm=(\overbrace{\ell+1,\dots,\ell+1}^{m}),\ 
\mu=(\overbrace{\ell,\dots,\ell}^{m-1},0)$ and $\peri=(m,-\ell)$.
\end{itemize}

\bigskip
\noindent
{\bf Acknowledgments.}
We thank H.~Tagawa for suggesting us a formula which leads a proof of Theorem \ref{thm1}.
We also thank K.~Nakada for discussion and valuable comments.

%%%%%%%%%%%%%%%%%%%%%%%%%%%%%%%%%%%%%%%%%%%%%%%%
\section{Cylindric diagrams}% 
%%%%%%%%%%%%%%%%%%%%%%%%%%%%%%%%%%%%%%%%%%%%%%%%%%%%%%%

For $\omega \in \mathbb{Z}_{\geqq 1} \times \mathbb{Z}_{\leqq -1}$,
we let $\mathbb{Z} \omega$ denote the subgroup of (the additive group) $\mathbb{Z}^2$ generated by $\omega$,
and define
$$
\mathcal{C}_\omega = \mathbb{Z}^2 / \mathbb{Z} \omega.
$$
We regard $\mathbb{Z}^2$ as a poset with the following partial order
$$
(a,b) \leq (a^\prime, b^\prime) \iff \hbox{$a \leqq a^\prime$ and $b \leqq b^\prime$ 
as integers}.
$$
Then, the cylinder $\C_\omega$ admits an induced poset structure, namely,
$$
x \le y
\iff \hbox{$\exists \tilde{x}, \tilde{y} \in \mathbb{Z}^2$ such that $\pi(\tilde{x})=x,\ \pi(\tilde{y})=y$ and
$\tilde{x} \le \tilde{y}$},
$$
where $\pi:\mathbb{Z}^2 \to \C_\omega$ is the natural projection.

\begin{df}
Let $(P,\leq)$ be a poset.
A subset $F$ of $P$ is called an \alert{order filter}
if the following condition holds:
$$x \in F,\ x \leq y \ \Rightarrow\ y \in F.$$
An order filter $F$ is said to be \alert{non-trivial} if $F\neq \emptyset$ nor $F\neq P$.
\end{df}

\begin{df}
Let  $\omega \in \mathbb{Z}_{\geqq1} \times \mathbb{Z}_{\leqq -1}$.
A non-trivial order filter $\cdiag$ of $\mathcal{C}_\omega$ is called a \alert{cylindric diagram}.
The inverse image $\pi^{-1}(\cdiag)\subset \Z^2$ %of a cylindric diagram
is called a \alert{periodic diagram} of period $\omega$.
An element of cylindric/periodic diagram is called a \alert{cell}.
\end{df}

Note that a cylindric diagram $\cdiag$ is  a poset and its order filter $\ofil$
is a cylindric diagram such thet $\ofil \subset \cdiag$.

%$\ell$-restricted skew diagram
%%%%%%%%%%%%%%%%%%%%%%%%%%%%%%%%%%%%%%%%%%%%%%%%%%%%%%%%%%%%%%%%%%%%%%%%%%%%%%
\begin{df}
Let $\theta\subset \mathcal{C}_\omega$ be a cylindric diagram and 
$\ofil$ an order filter of $\theta$.
The set-difference $$\theta/\eta:=\theta\backslash\eta$$  
is called a \alert{cylindric skew diagram}.
The inverse image $\pi^{-1}(\cdiag/\ofil)\subset \Z^2$ is called a 
\alert{periodic skew diagram} of period $\omega$.
\end{df}

We sometimes parameterize periodic/cylindric diagrams by $\ell$-restricted partitions:
%
%%%%%%%%%%%%%%%%%%%%%%%%%%%%%%%%%%%%%%%%%%%%%
\begin{df}
Let $m,\ell\in\Z_{\geqq1}$.
An integer sequence $\lambda=(\lambda_1,\dots,\lambda_m)$ is an \alert{$\ell$-restricted generalized partition of length $m$}
if it satisfies the following conditions:
$$
\lambda_1 \geqq \cdots \geqq \lambda_\ell,\quad
\lambda_1-\lambda_\ell \leqq \ell.
$$
We 
denote by $\P_{m,\ell}$ the set of $\ell$-restricted 
generalized partitions of length $m$.
(Note that we allow $\lambda_i$ to be negative).
\end{df}

For a partition $\lm=(\lm_1,\dots,\lm_m)$, we denote by the same symbol $\lm$
the following subset of $\Z^2$:
$$\lm=\{(a,b)\in\Z^2\mid 1\leqq a\leqq m,\ 1\leqq b\leqq \lm_a\},$$
which is called a (Young) diagram associated with $\lm$.

Let $\lambda=(\lambda_1,\dots,\lambda_m) \in \P_{m,\ell}$.
We define
\begin{align*}
\semid\lambda &= \{ (a,b) \in \mathbb{Z}^2 \mid 1 \leqq a \leqq m,\ b \leqq \lambda_a \}, \\
\perid{\lambda} &=\perid{\lm}_{(m,-\ell)}= \semid\lambda + \mathbb{Z} (m,-\ell),\\
\cyl\lm &=\cyl\lm_{(m,-\ell)}=\pi(\perid\lm).
\end{align*}
Note that  
$\semid\lambda = \perid{\lambda} \cap ([1,m] \times \mathbb{Z})$ 
and $\lm$ is a fundamental domain of $\perid\lambda$ with respect to the action of
$\mathbb{Z} (m,-\ell)$.

It is easy to see that 
$\perid{\lambda}$ is a periodic diagram of period  $(m,-\ell)$
and $\cyl\lm$ is a cylindric diagram.
Moreover,
 any periodic (resp. cylindric)
  diagram of period $(m,-\ell)$ is of the form $\perid\lambda$
  (resp. $\pi(\perid{\lambda})$ for some 
 $\lambda \in \P_{m,\ell}$.

%%%%%%%%%%%%%%%%%%%%%%%%%%%%%%%%%%%%%%%%%%%%%%%%%%%%%%%%%%%%%%%%%%%%%%%%%%%%%%%%

\begin{figure}[h]
\begin{small}
\begin{center}
\begin{tikzpicture}[scale=0.4, rotate=270]
\fill[pink] (0,-5) |- (2,2) |- (3,0) |- (4,-1) |- (0,-5) -- cycle;
\fill[yellow] (0,2) |- (1,4) |- (3,3) |- (4,1) |- (3,-1) |- (2,0) |- (1,2) -- (0,2) -- cycle; 
\foreach \i/\l in {-3/8, -2/7, -1/7, 0/5, 1/4, 2/3, 3/3, 4/1, 5/0, 6/-1, 7/-1, 8/-3}
\draw (\i,-5) grid (\i-1,\l);
\draw[ultra thick] (-4,-5) |- (-3,8) |- (-1,7) |- (0,5) -- (0,-5);
\draw[ultra thick] (0,-5) |- (1,4) |- (3,3) |- (4,1) -- (4,-5);
\draw[ultra thick] (4,-5) |- (5,0) |- (7,-1) |- (8,-3) -- (8,-5);

\draw[dotted, ultra thick] (-5,-4) -- +(-1,0);
\draw[dotted, ultra thick] (-5,0) -- +(-1,0);
\draw[dotted, ultra thick] (-5,4) -- +(-1,0);

\draw[dotted, ultra thick] (-2,-6) -- +(0,-1);
\draw[dotted, ultra thick] (2,-6) -- +(0,-1);
\draw[dotted, ultra thick] (6,-6) -- +(0,-1);

\draw[dotted, ultra thick] (9,-4) -- +(1,0);

\draw[->, thick] (2.5,4) -- +(4,-4);
\node[anchor=west] at (4.5,2) {$(4,-4)$};

\end{tikzpicture}
\end{center}
\end{small}
\caption{The semi-infinite periodic diagram of period $(4,-4)$ associate with $\lambda=(5,4,4,2)$.}\label{fig1}
\end{figure}

%%%%%%%%%%%%%%%%%%%%%%%%%%%%%%%%%%%%%%
\section{Linear extensions}
%%%%%%%%%%%%%%%%%%%%%%%%%%%%%%%%%%%%%%%

For two integers $a,b$,
we  use the following notation:
$$[a,b]=\{x\in \Z\mid a\leqq x\leqq b\}.$$
%%%%%%%%%%%%%%%%%%%%%%%%%%%%%%%%%%%%%%%%%%%%%%%%%%%%%%%%%%%%%%%%%%%%%%%%
\begin{df}
For a poset $P$ such that $|P|=n$,
 a \alert{linear extension} (or a  \alert{reverse standard tableau})
of $P$ is a bijection 
$\lex:P \to [1,n]$ satisfying
$$
x<y \Longrightarrow \lex(x)<\lex(y).
$$
Let $\LE(P)$ denote the set of all linear extensions of 
$P$.
\end{df}
%%%%%%%%%%%%%%%%%%%%%%%%%%%%%%%%%%%%%%%%%%%%%%%%%%%%%%%%%%%%%%%%
%
Let $\lm,\mu$ be partitions such that $\lm\supset \mu$ and $|\lm/\mu|=n$.

It is easy to see that 
a bijection $\lex:\lm/\mu\to [1,n]$ is a linear extension on the finite skew diagram $\lm/\mu$
if and only if the following conditions hold:

\begin{enumerate}[(1)]
\item $\lex(a,b)>\lex(a,b+1)$ whenever $(a,b),(a,b+1)\in \lm/\mu$.

\item $\lex(a,b)>\lex(a+1,b)$ whenever $(a,b),(a+1,b)\in \lm/\mu$.
\end{enumerate}

Fix $m,\ell\in \Z_{\geqq1}$.
%%%%%%%%%%%%%%%%%%%%%%%%%%%%%%%%%%%%%
\begin{df}
Let $\lm,\mu\in \P_{m,\ell}$.
A linear extension $\lex$ on $\lm/\mu$ is \alert{$\ell$-restricted}
if it satisfies
$$\lex(1,b)<\lex(m,b-\ell)\ \text{
whenever }(1,b),(m,b-\ell)\in \lm/\mu.$$

We denote by  $\LE_\ell(\lm/\mu)$ the set of all $\ell$-restricted 
linear extensions of $\lm/\mu$.
\end{df}
%%%%%%%%%%%%%%%%%%%%%%%%%%%%%%%%%%%%%%%

Note that the projection $\pi:\Z^2\to\C_\peri$ gives a bijection
$\lm/\mu\to \pi(\lm/\mu)=\cy\lm/\cy\mu$ and
that a map $\lex: \lm/\mu\to [1,n]$ induces a map
$\lex: \cy\lm/\cy\mu\to [1,n]$.
%%%%%%%%%%%%%%%%%%%%%%%%%%%%%%%%%%%%%%%%%%%%%%%%
\begin{lem}
Let $\lm,\mu\in \P_{m,\ell}$.
A bijection $\lex:\lm/\mu\to [1,n]$
induces a linear extension on the cylindric skew diagram $\cy\lm/\cy\mu$ 
if and only if $\lex\in \LE_\ell(\lm/\mu)$.
Namely, the set $\LE(\cy\lm/\cy\mu)$ and $\LE_\ell(\lm/\mu)$
are in one to  one correspondence.
\end{lem}

%%%%%%%%%%%%%%%%%%%%%%%%%%%%%%%%%%%%%%%%%%%%%%%%%%%%%
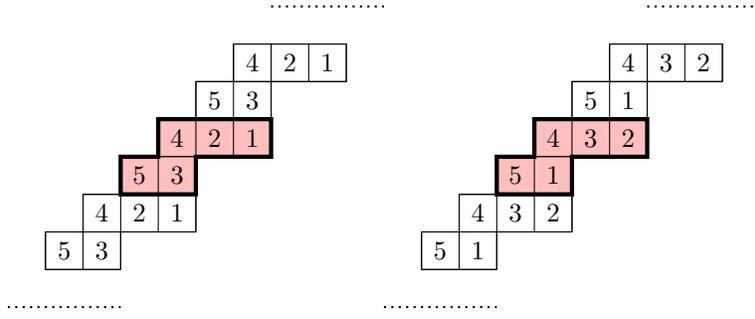
\begin{figure}[h]
\begin{center}
\begin{tikzpicture}[rotate=270, scale=0.5]
\foreach \n/\l/\m in {1/4/1, 2/2/0}{
\fill[pink] (\n,\m) rectangle (\n-1,\l);
\draw (\n,\m) grid (\n-1,\l);
}
\draw[ultra thick] (0,1) |- (1,4) |- (2,2) |- (1,0) |- (0,1) -- cycle;

\begin{scope}[xshift=-0.5cm, yshift=-0.5cm]
\foreach \i/\j/\t in {1/2/4, 1/3/2, 1/4/1, 2/1/5, 2/2/3}
\node at (\i,\j) {$\t$};
\end{scope}

\begin{scope}[xshift=2cm, yshift=-2cm]
\foreach \n/\l/\m in {1/4/1, 2/2/0}
\draw (\n,\m) grid (\n-1,\l);
\draw[thick, dotted]  (3,-1) -- +(0,3);
\draw (0,1) |- (1,4) |- (2,2) |- (1,0) |- (0,1) -- cycle;

\begin{scope}[xshift=-0.5cm, yshift=-0.5cm]
\foreach \i/\j/\t in {1/2/4, 1/3/2, 1/4/1, 2/1/5, 2/2/3}
\node at (\i,\j) {$\t$};
\end{scope}
\end{scope}

\begin{scope}[xshift=-2cm, yshift=2cm]
\foreach \n/\l/\m in {1/4/1, 2/2/0}
\draw (\n,\m) grid (\n-1,\l);
\draw[thick, dotted]  (-1,2) -- +(0,3);
\draw (0,1) |- (1,4) |- (2,2) |- (1,0) |- (0,1) -- cycle;

\begin{scope}[xshift=-0.5cm, yshift=-0.5cm]
\foreach \i/\j/\t in {1/2/4, 1/3/2, 1/4/1, 2/1/5, 2/2/3}
\node at (\i,\j) {$\t$};
\end{scope}
\end{scope}

\begin{scope}[yshift=10cm]
\foreach \n/\l/\m in {1/4/1, 2/2/0}{
\fill[pink] (\n,\m) rectangle (\n-1,\l);
\draw (\n,\m) grid (\n-1,\l);
}
\draw[ultra thick] (0,1) |- (1,4) |- (2,2) |- (1,0) |- (0,1) -- cycle;

\begin{scope}[xshift=-0.5cm, yshift=-0.5cm]
\foreach \i/\j/\t in {1/2/4, 1/3/3, 1/4/2, 2/1/5, 2/2/1}
\node at (\i,\j) {$\t$};
\end{scope}

\begin{scope}[xshift=2cm, yshift=-2cm]
\foreach \n/\l/\m in {1/4/1, 2/2/0}
\draw (\n,\m) grid (\n-1,\l);
\draw[thick, dotted]  (3,-1) -- +(0,3);
\draw (0,1) |- (1,4) |- (2,2) |- (1,0) |- (0,1) -- cycle;

\begin{scope}[xshift=-0.5cm, yshift=-0.5cm]
\foreach \i/\j/\t in {1/2/4, 1/3/3, 1/4/2, 2/1/5, 2/2/1}
\node at (\i,\j) {$\t$};
\end{scope}
\end{scope}

\begin{scope}[xshift=-2cm, yshift=2cm]
\foreach \n/\l/\m in {1/4/1, 2/2/0}
\draw (\n,\m) grid (\n-1,\l);
\draw[thick, dotted]  (-1,2) -- +(0,3); % {$\vdots$};
\draw (0,1) |- (1,4) |- (2,2) |- (1,0) |- (0,1) -- cycle;

\begin{scope}[xshift=-0.5cm, yshift=-0.5cm]
\foreach \i/\j/\t in {1/2/4, 1/3/3, 1/4/2, 2/1/5, 2/2/1}
\node at (\i,\j) {$\t$};
\end{scope}
\end{scope}
\end{scope}
\end{tikzpicture}
\end{center}
\caption{Two images of the map from 
$\cy\lm/\cy\mu$ to $\{1,2,3,4,5\}$ with $\lm=(4,2)$, $\mu=(1,0)$ and $\ell=2$.
The left is a linear extension, but the right is not a linear extension.}
\end{figure}

%%%%%%%%%%%%%%%%%%%%%%%%%%%%%%%%%%%%%%%%%%%%%%%%%%%%%%%
\section{Excited diagrams}
%%%%%%%%%%%%%%%%%%%%%%%%%%%%%%%%%%%%%%%%%%%%%%%%%%%%%%%%%%

In this section, we fix $\ell,m \in \mathbb{Z}_{\geqq 1}$.

\begin{df}
Let $\lambda$ and $\mu $ be two partitions such that $\lambda \supset \mu$.

\begin{enumerate}[(1)]
\item 
Let $D$ be a subset of $\lambda$.
A cell $y=(a,b) \in D$ is said to be \alert{$D$-active} 
if 
$$(a+1,b),\ (a,b+1),\ (a+1,b+1)\in \lambda\setminus D.$$

\item
For a $D$-active cell $y=(a,b)$, we put
$$
D_y = D \setminus \{y\} \cup \{x\},
$$
where $x=(a+1,b+1)$.
The replacement from $D$ to $D_y$ is called an \alert{elementary excitation} at $y$.

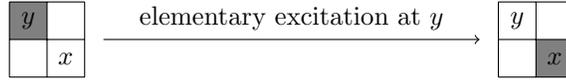
\begin{center}
\begin{tikzpicture}[scale=0.5]
\fill[gray] (0,1) rectangle +(1,1);
\draw (0,0) grid (2,2);
\node at (0.5,1.5) {$y$};
\node at (1.5,0.5) {$x$};

\draw[->] (2.5,1) -- +(10,0);
\node[above] at (7.5,1) {elementary excitation at $y$};

\begin{scope}[xshift=13cm]
\fill[gray] (1,0) rectangle +(1,1);
\draw (0,0) grid (2,2);
\node at (0.5,1.5) {$y$};
\node at (1.5,0.5) {$x$};
\end{scope}

\end{tikzpicture}
\end{center}

\item An \alert{excited diagram} of $\mu$ in $\lambda$ is a subset of $\lambda$ obtained from $\mu$ after a sequence of elementary excitations on active cells.
Let $\E_\lambda(\mu)$ denote the set of all excited diagrams of $\mu$ in $\lambda$.
\end{enumerate}
\end{df}

%%%%%%%%%%%%%%%%%%%%%%%%%%%%%%%%%%%%%%%%%%

We extend the concept of excited diagrams to the case of cylindric/periodic diagrams.

\begin{df}
Let $\lambda$ and $\mu$ be two $\ell$-restricted partitions of length $m$
such that $\lambda \supset \mu$.
Put $\peri=(m,-\ell)$.

\begin{enumerate}[(1)]
\item
Let $D$ be a ``periodic'' subset of $\perid{\lambda}$
(i.e., $D+\peri=D$).
A cell $(a,b) \in D$ is \alert{$D$-active}
if
$$
(a+1,b), (a,b+1), (a+1,b+1) \in \perid{\lambda} \setminus D.
$$

\item
Let $(a,b) \in D$ be a $D$-active cell.
Put
$$
D_{(a,b)} :=
D \setminus ((a,b)+\Z\peri) \cup ((a+1,b+1)+\Z\peri).
$$
The replacement from $D$ to $D_{(a,b)}$ is called a \alert{periodic elementary excitation}.

\item
A \alert{periodic excited diagram} of $\perid{\mu}$ in $\perid{\lambda}$ is
a subset of $\perid{\lambda}$ obtained from $\perid{\mu}$ after a sequence of periodic elementary excitations on active cells,
and the whole set is denoted by $\E_{\perid{\lambda}}(\perid{\mu})$.

\item
For $D \in \E_{\perid{\lambda}}(\perid{\mu})$,
$\pi(D)$ is called a \alert{cylindric excited diagram},
and the whole set is denoted by $\E_{\cy\lm}(\cy\mu)$:
$$
\E_{\cy\lm}(\cy\mu) = \{ \pi(D) \mid D \in \E_{\perid\lm}(\perid\mu)\}.
$$
\end{enumerate}
\end{df}

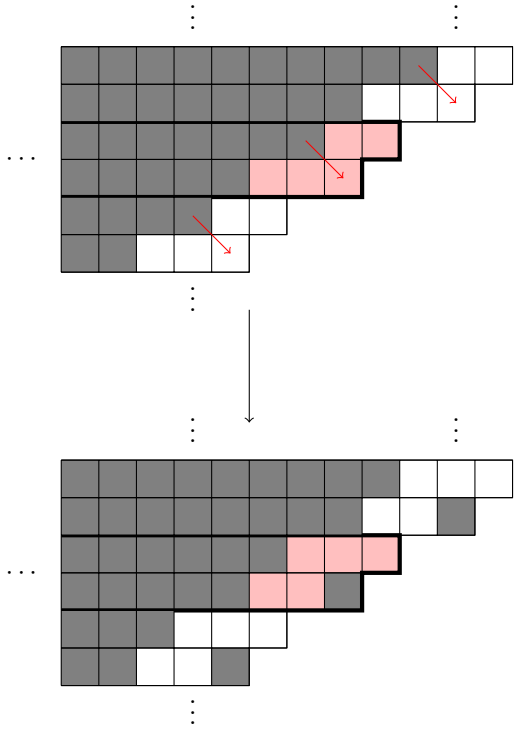
\begin{figure}[h]
\begin{center}
\begin{tikzpicture}[rotate=270,scale=0.5]
\fill[gray] (1,-5) rectangle (0,2);
\fill[gray] (2,-5) rectangle (1,0);
\fill[pink] (1,2) rectangle (0,4);
\fill[pink] (2,0) rectangle (1,3);
\draw (1,-5) grid (0,4);
\draw (2,-5) grid (1,3);
\draw[ultra thick] (2,-5) |- (1,3) |- (0,4) -- (0,-5);
\draw[red,->] (0.5,1.5) -- +(1,1);

\begin{scope}[xshift=2cm,yshift=-3cm]
\fill[gray] (1,-2) rectangle (0,2);
\fill[gray] (2,-2) rectangle (1,0);
\draw (1,-2) grid (0,4);
\draw (2,-2) grid (1,3);
\draw[] (2,-2) |- (1,3) |- (0,4) -- (0,-2);
\draw[red,->] (0.5,1.5) -- +(1,1);
\node at (2.5,1.5) {$\vdots$};
\end{scope}

\begin{scope}[xshift=-2cm,yshift=3cm]
\fill[gray] (1,-8) rectangle (0,2);
\fill[gray] (2,-8) rectangle (1,0);
\draw (1,-8) grid (0,4);
\draw (2,-8) grid (1,3);
\draw[] (2,-8) |- (1,3) |- (0,4) -- (0,-8);
\draw[red,->] (0.5,1.5) -- +(1,1);
\node at (-1,2.5) {$\vdots$};
\node at (-1,-4.5) {$\vdots$};
\end{scope}

\node at (1,-6) {$\cdots$};
\draw[->] (5,0) -- +(3,0);

\begin{scope}[xshift=11cm]
\fill[gray] (1,-5) rectangle (0,1);
\fill[gray] (2,-5) rectangle (1,0);
\fill[gray] (1,2) rectangle +(1,1);
\fill[pink] (1,1) rectangle (0,4);
\fill[pink] (2,0) rectangle (1,3);
\fill[gray] (1,2) rectangle +(1,1);
\draw (1,-5) grid (0,4);
\draw (2,-5) grid (1,3);
\draw[ultra thick] (2,-5) |- (1,3) |- (0,4) -- (0,-5);

\begin{scope}[xshift=2cm,yshift=-3cm]
\fill[gray] (1,-2) rectangle (0,1);
\fill[gray] (2,-2) rectangle (1,0);
\fill[gray] (1,2) rectangle +(1,1);
\draw (1,-2) grid (0,4);
\draw (2,-2) grid (1,3);
\draw[] (2,-2) |- (1,3) |- (0,4) -- (0,-2);
\node at (2.5,1.5) {$\vdots$};
\end{scope}

\begin{scope}[xshift=-2cm,yshift=3cm]
\fill[gray] (1,-8) rectangle (0,1);
\fill[gray] (2,-8) rectangle (1,0);
\fill[gray] (1,2) rectangle +(1,1);
\draw (1,-8) grid (0,4);
\draw (2,-8) grid (1,3);
\draw[] (2,-8) |- (1,3) |- (0,4) -- (0,-8);
\node at (-1,2.5) {$\vdots$};
\node at (-1,-4.5) {$\vdots$};
\end{scope}

\node at (1,-6) {$\cdots$};

\end{scope}

\end{tikzpicture}
\end{center}
\caption{A periodic elementary excitation.}
\end{figure}

\begin{rem}
It is easy to see that
if $u \in D$ is $D$-active,
then $u+k\peri$ is also $D$-active for any $k \in \Z$,
and that $D_u$ is periodic.
\end{rem}

%%%%%%%%%%%%%%%%%%%%%%%%%%%%%%%%%%%%%%%%%%%%%%%%%%%%%%%%%%%%%%%%%%%%%%%%%%%%%%%%
\section{Conjectural hook formula for cylindric skew diagrams}
%%%%%%%%%%%%%%%%%%%%%%%%%%%%%%%%%%%%%%%%%%%%%%%%%%%%%%%%%%%%%%%%%%%%%

\begin{df}
Let $\lm$ be a partition.
For a cell $x$ of the corresponding
finite Young diagram $\lm$,
the \alert{hook} $H_\lm(x)$ of $x$ in $\lm$ is given by
$$
H_\lm(x) = \lm \cap \Bigl(\{ x+(k,0) \mid k \in \Z_{\geqq0} \} \cup \{ x+(0,k) \mid k \in \Z_{\geqq1} \}\Bigr).
$$
and the \alert{hook length} $h_\lm(x)$ is the number of cells of $H_\lm(x)$.
\end{df}

Fix $m,\ell \in \Z_{\geqq1}$ and put $\peri=(m,-\ell)$.
%%%%%%%%%%%%%%%%%%%%%%%%%%%%%%%%%%%%%%%%%%
\begin{df}
Let $\lm \in \P_{m,\ell}$.
For a cell $x \in \perid\lm$,
define the hook $H_{\perid\lm}(x)$ of $x$ in $\perid\lm$ by
$$
H_{\perid{\lambda}}(x) =
\perid\lm \cap \Bigl(\{ x+(k,0) \mid k \in \Z_{\geqq0} \} \cup \{ x+(0,k) \mid k \in \Z_{\geqq1} \}\Bigr).
$$
The number $h_{\perid{\lambda}}(x)$ of cells of $H_{\perid{\lambda}}(x)$ is called the \alert{hook length} in $\perid\lm$. 
\end{df}
%%%%%%%%%%%%%%%%%%%%%%%%%%%%%%%%%%

%%%%%%%%%%%%%%%%%%%%%%%%%%%%%%%%%%%%
\begin{df}\label{df;chook}
For a cell $x \in \cy{\lambda}$,
define 
$$
h_{\cyl{\lambda}}(x) := h_{\perid{\lambda}}(y),
$$
where $y \in \pi^{-1}(x)$.
\end{df}
%%%%%%%%%%%%%%%%%%%%%%%%%%%%%%%%%%%%%%%%%%%%%%%%%%%%%%%%%%%%%
Note that $h_{\cyl{\lambda}}(x)$ is well-defined
since
$$
h_{\perid{\lambda}}(y+\peri) = h_{\perid{\lambda}}(y).
$$

\begin{figure}[h]
	\begin{center}
		\begin{tikzpicture}[rotate=270,scale=0.5]
		%\foreach \i in {-4,...,8}
		%\node at (-0.5,\i-3.5) {\i};
		
		\foreach \n/\l in {1/5, 2/3, 3/2}
		\draw (\n,-8) grid (\n-1,\l);
		\draw[] (0,-8) |- (1,5) |- (2,3) |- (3,2) -- (3,-8);
		\node at (1.5,-9) {$\cdots$};
		\node at (-1,0) {$\vdots$};
		
		\begin{scope}[xshift=3cm,yshift=-3cm]
		\foreach \n/\l in {1/5, 2/3, 3/2}
		\fill[pink] (\n,-5) rectangle (\n-1,\l);
		
		\fill[gray] (2,0) rectangle (1,3);
		\fill[gray] (4,0) rectangle (2,1);
		\node at (1.5,0.5) {$x$};
		
		\foreach \n/\l in {1/5, 2/3, 3/2}
		\draw (\n,-5) grid (\n-1,\l);
		\draw[ultra thick] (0,-5) |- (1,5) |- (2,3) |- (3,2) -- (3,-5);
		\node at (1.5,-6) {$\cdots$};
		\end{scope}
		
		\begin{scope}[xshift=6cm,yshift=-6cm]
		\foreach \n/\l in {1/5, 2/3, 3/2}
		\draw (\n,-2) grid (\n-1,\l);
		\draw[] (0,-2) |- (1,5) |- (2,3) |- (3,2) -- (3,-2);
		\node at (1.5,-3) {$\cdots$};
		\node at (4,0) {$\vdots$};
		\node at (-7,0) {$\vdots$};
		\end{scope}
		\end{tikzpicture}
	\end{center}
\caption{
The hook $H_{\perid{\lm}}(x)$ of $x$ in the periodic diagram $\perid\lm$.
The hook length $h_{\perid\lm}(x)=5$.
}
\label{cylper}
\end{figure}

\if0
For $\lm,\mu \in \P_{m,\ell}$ with $\lm \supset \mu$,
we define $f^{\lm/\mu}$ (resp. $f^{\cyl\lm/\cyl\mu}$)
as the number of linear extensions of $\lm/\mu$ (resp. $\cyl\lm/\cyl\mu$):
\begin{align*}
f^{\lm/\mu} &= |\LE(\lm/\mu)|,\\
f^{\cyl\lm/\cyl\mu} &= |\LE(\cyl\lm/\cyl\mu)| = |\LE_\ell(\lm/\mu)|.
\end{align*}
\fi

For a skew diagram $\lm/\mu$,
we denote by $f^{\lm/\mu}$ the number of linear extensions of $\lm/\mu$:
$$
f^{\lm/\mu} = |\LE(\lm/\mu)|.
$$

%%%%%%%%%%%%%%%%%%%%%%%%%%%%%%%%
\begin{thm}[\cite{Nar}]\label{thm;NHF}
Let $\lambda$ and $\mu$ be partitions with $\lambda \supset \mu$ and $|\lambda/\mu|=n$.
Then, %the number $f^{\lambda/\mu}$ of linear extensions of the skew Young diagram $\lambda/\mu$ is given by
\begin{equation}\label{eq;NHF}
f^{\lambda/\mu} = n! \sum_{D \in \mathcal{E}_\lambda(\mu)} \prod_{x \in \lambda \setminus D} \frac{1}{h_\lambda(x)},
\end{equation}
where $\mathcal{E}_\lambda(\mu)$ is the set of all excited diagrams of $\mu$ in $\lambda$,
and $h_\lm(x)$ is the hook length of $x$ in $\lm$.
\end{thm}
%%%%%%%%%%%%%%%%%%%%%%%%%%%%%

For a cylilndric skew diagram $\cyl\lm/\cyl\mu$ of period $(m,-\ell)$,
we denote by $f^{\cyl\lm/\cyl\mu}$ the number of linear extensions of $\cyl\lm/\cyl\mu$:
$$
f^{\cyl\lm/\cyl\mu} = |\LE(\cyl\lm/\cyl\mu)| = |\LE_\ell(\lm/\mu)|.
$$

%%%%%%%%%%%%%%%%%%%%%%%%%%%%%%%%%%%%%%%%%%%%%%%%%%%%%%%%%%%%%%%%%%%%%%%%%
\begin{conj}\label{conj2}
Let $m,\ell \in \Z_{\geqq1}$
and $\lm,\mu \in \P_{m,\ell}$ such that $\lm \supset \mu$.
Put $n = |\lm/\mu| = |\cy\lm / \cy\mu|$.
Then, %the number $f^{\cy\lm / \cy\mu}$ of linear extensions of $\cy\lm/\cy\mu$ is given by
\begin{equation}\label{eq;conj2}
f^{\cy\lm/\cy\mu} = n! \sum_{D \in \mathcal{E}_{\cy\lm}(\cy\mu)} \prod_{x \in \cy\lm \setminus D}
\frac{1}{h_{\cy\lm}(x)},
\end{equation}
where $\mathcal{E}_{\cy\lm}(\cy\mu)$ is the set of all cylindric excited diagrams of $\cy\mu$ in $\cy\lm$,
and $h_{\cy\lm}(x)$ is the hook length of $x$ in $\cy\lm$.
\end{conj}
%%%%%%%%%%%%%%%%%%%%%%%%%%%%%%%%%%%%%%%%%%%%%%%%%%%%%%%%%%%%%%%%%%%%%%%%%%%%

%
\begin{rem}
If $\ell\geqq\lm_1$,
then
\begin{gather*}
\E_\cyl\lm(\cyl\mu) = \E_\lm(\mu), \\
h_\cyl\lm(\pi(x)) = h_\lm(x),
\end{gather*}
and hence
Conjecture \ref{conj2} follows from Theorem \ref{thm;NHF}.
\fin
\end{rem}

\begin{ex}
Let us see the simplest non-trivial example.
Let $\lambda=(2)$,  $\mu=(0) \in \P_{1,1}$.
Then $\cyl{\lambda}/\cyl{\mu}$ has just one linear extension.
%%%%%%%%%%%%%%%%%%%%%%%%%%%%%%%%%%%%%%%%%%%%%%%%%%%%%%%%%%%%%%%%%%%%%%%%%%%%%%%%

The hook length on fundamental domain of $\cyl{\lambda}$ are as follows:
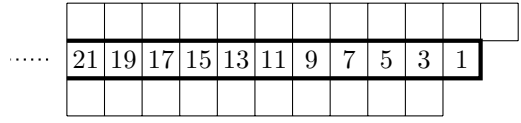
\begin{center}
\begin{tikzpicture}[scale=0.5]
\draw (-9,1) grid (3,2);
\draw (-9,0) grid (2,1);	\draw[ultra thick] (-9,0) -| (2,1) -- (-9,1);
\draw (-9,-1) grid (1,0);
\draw[thick, dotted] (-9.5,0.5) -- +(-1,0);

\foreach \x/\h in {2/1, 1/3, 0/5, -1/7, -2/9, -3/11, -4/13, -5/15, -6/17, -7/19, -8/21}
\node at (\x-0.5,0.5) {$\h$};
\end{tikzpicture}
\end{center}
The excited diagrams of $\cyl{\mu}$ in $\cyl{\lambda}$ are as follows:
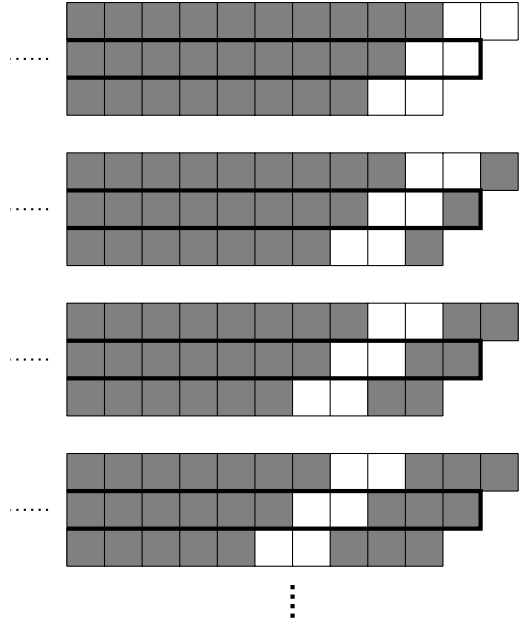
\begin{center}
\begin{tikzpicture}[scale=0.5]
\fill[gray] (-9,1) rectangle (1,2);
\fill[gray] (-9,0) rectangle (0,1);
\fill[gray] (-9,-1) rectangle (-1,0);
\draw (-9,1) grid (3,2);
\draw (-9,0) grid (2,1);	\draw[ultra thick] (-9,0) -| (2,1) -- (-9,1);
\draw (-9,-1) grid (1,0);
\draw[thick, dotted] (-9.5,0.5) -- +(-1,0);

\begin{scope}[yshift=-4cm]
\fill[gray] (-9,1) rectangle (0,2);		\fill[gray] (2,1) rectangle (3,2);
\fill[gray] (-9,0) rectangle (-1,1);	\fill[gray] (1,0) rectangle (2,1);
\fill[gray] (-9,-1) rectangle (-2,0);	\fill[gray] (0,-1) rectangle (1,0);
\draw (-9,1) grid (3,2);
\draw (-9,0) grid (2,1);	\draw[ultra thick] (-9,0) -| (2,1) -- (-9,1);
\draw (-9,-1) grid (1,0);
\draw[thick, dotted] (-9.5,0.5) -- +(-1,0);
\end{scope}

\begin{scope}[yshift=-8cm]
\fill[gray] (-9,1) rectangle (-1,2);		\fill[gray] (1,1) rectangle (3,2);
\fill[gray] (-9,0) rectangle (-2,1);	\fill[gray] (0,0) rectangle (2,1);
\fill[gray] (-9,-1) rectangle (-3,0);	\fill[gray] (-1,-1) rectangle (1,0);
\draw (-9,1) grid (3,2);
\draw (-9,0) grid (2,1);	\draw[ultra thick] (-9,0) -| (2,1) -- (-9,1);
\draw (-9,-1) grid (1,0);
\draw[thick, dotted] (-9.5,0.5) -- +(-1,0);
\end{scope}

\begin{scope}[yshift=-12cm]
\fill[gray] (-9,1) rectangle (-2,2);		\fill[gray] (0,1) rectangle (3,2);
\fill[gray] (-9,0) rectangle (-3,1);	\fill[gray] (-1,0) rectangle (2,1);
\fill[gray] (-9,-1) rectangle (-4,0);	\fill[gray] (-2,-1) rectangle (1,0);
\draw (-9,1) grid (3,2);
\draw (-9,0) grid (2,1);	\draw[ultra thick] (-9,0) -| (2,1) -- (-9,1);
\draw (-9,-1) grid (1,0);
\draw[thick, dotted] (-9.5,0.5) -- +(-1,0);
\draw[ultra thick, dotted] (-3,-1.5) -- +(0,-1);
\end{scope}
\end{tikzpicture}
\end{center}
Therefore, by computing the right hand side of (\ref{eq;conj2}),
\begin{align*}
f^{\cyl{\lambda}/\cyl{\mu}}
&= 2! \left(\frac{1}{1 \cdot 3} + \frac{1}{3 \cdot 5} + \frac{1}{5 \cdot 7} + \frac{1}{7 \cdot 9} + \cdots \right) \\
&= 2! \cdot \sum_{k=0}^\infty \frac{1}{(2k+1)(2k+3)} \\
&= 2! \cdot \frac{1}{2} \cdot \sum_{k=0}^\infty 
\left( \frac{1}{2k+1} - \frac{1}{2k+3} \right) 
=1.
\end{align*}
%%%%%%%%%%%%%%%%%%%%%%%%%%%%%%%%%%%%%%%%%%%%%%%%%%%%%%%%%%%%%%%%%%%%%%%%%%%%%%%%
%%%%%%%%%%%%%%%%%%%%%%%%%%%%%%%%%%%%%%%%%%%%%%%%%%%%%%%%%%%%%%%%%%%%%%%%%%%%%%%%
\fin
\end{ex}

By similar case by case compilation, 
we have confirmed Conjecture \ref{conj2} for any shape with $n \leqq 4$.

In the rest, we denote the right hand side of \eqref{eq;conj2} by $g^{\cyl\lm/\cyl\mu}$:
$$
g^{\cyl\lm/\cyl\mu} = n! \sum_{D \in \mathcal{E}_{\cy\lm}(\cy\mu)} \prod_{x \in \cy\lm \setminus D}
\frac{1}{h_{\cy\lm}(x)}.
$$
The proofs for the following two theorems will be given in the later sections:
%%%%%%%%%%%%%%%%%%%%%%%%%%%%%%%%%%%%%%%%%%%%%%%%%%%%%%%%%%%%%%%%%%%%%%
\begin{thm}\label{thm1}
{\rm{(Bar cases)}}\ \
Let $n, \ell\in\Z_{\geqq1}$.
Put $\lm=(n)$ and $\mu=(0)$, which belong to $\P_{1,\ell}$,
and let $\cyl\lm=\cyl\lm_{(1,-\ell)}$ and $\cyl\mu=\cyl\mu_{(1,-\ell)}$ be 
the corresponding cylindric diagrams.
Then
\begin{equation}\label{g=f1}
f^{\cyl\lm/\cyl\mu} =1=g^{\cyl\lm/\cyl\mu} 
\end{equation}
\end{thm}
%%%%%%%%%%%%%%%%%%%%%%%%%%%%%%%%%%%%%%%%%%%%%%%%%%%%%%%%%%%
\begin{thm}\label{thm2} {\rm{(Hook cases)}}\ \ 
Let $\ell, m\in\Z_{\geqq1}$.
Put $\lambda=((\ell+1)^m)$ and $\mu=(\ell^{m-1},0)$, 
which belong to $\P_{m,\ell}$,
and let $\cyl\lm=\cyl\lm_{(m,-\ell)}$ and $\cyl\mu=\cyl\mu_{(m,-\ell)}$  
be the corresponding cylindric diagrams.
Then
\begin{equation}\label{g=f2}
f^{\cyl\lm/\cyl\mu} =\binom{\ell+m-2}{m-1}=g^{\cyl\lm/\cyl\mu} 
\end{equation}
\end{thm}

%%%%%%%%%%%%%%%%%%%%%%%%%%%%%%%%%%%%%%%%%%%%%%%%%%%%%%%%%%%%%%%%%%%%%%%%%%%%%%%%%%%%%%%%%%%%
%Theorem \ref{thm1} (1) is proved just by case by case computation.

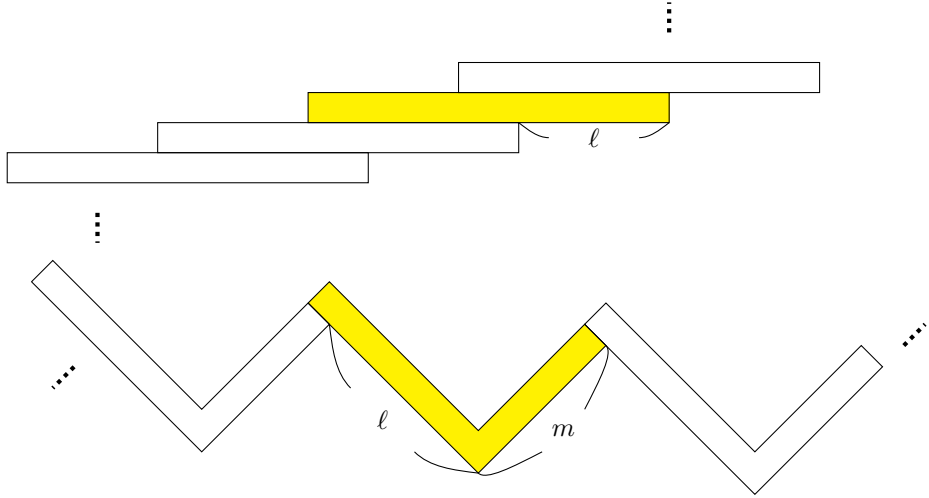
\begin{figure}[h]
\begin{center}
\begin{tikzpicture}[scale=0.4]

\fill[yellow] (0,0) rectangle +(12,1);

\draw[ultra thick, dotted] (12,3) -- +(0,1);
\draw (5,1) rectangle +(12,1);
\draw (0,0) rectangle +(12,1);
\draw (7,0) sin +(1,-0.5);
\draw (12,0) sin +(-1,-0.5);
\draw (-5,-1) rectangle +(12,1);
\draw (-10,-2) rectangle +(12,1);
\node at (9.5,-0.5) {$\ell$};
\draw[ultra thick, dotted] (-7,-3) -- +(0,-1);

\begin{scope}[yshift=-6cm, rotate=315]
\fill[yellow] (0,0) -| (8,6) -| (7,1) -| (0,0) -- cycle;
\draw (0,0) -| (8,6) -| (7,1) -| (0,0);
\node at (9,3) {$m$};
\draw (8,0) cos +(1,2);	\draw (8,6) cos +(1,-2);
\node at (4.5,-1) {$\ell$};
\draw (1,0) sin +(2,-1);	\draw (8,0) sin +(-2,-1);
\begin{scope}[xshift=-7cm,yshift=-6cm]
\draw (0,0) -| (8,6) -| (7,1) -| (0,0);
\draw[ultra thick, dotted] (3,-1) -- +(0,-1);
\end{scope}
\begin{scope}[xshift=7cm,yshift=6cm]
\draw (0,0) -| (8,6) -| (7,1) -| (0,0);
\draw[ultra thick, dotted] (8,7) -- +(0,1);
\end{scope}
\end{scope}
\end{tikzpicture}
\end{center}
\caption{The shapes indicated by Theorem \ref{thm1} and \ref{thm2}, respectively.}\label{fig601}
\end{figure}

It is easy to see the following:
%%%%%%%%%%%%%%%%%%%%%%%%%%%%%
\begin{prop}\label{prop:shift}
Let $m,\ell\in \Z_{\geqq1}$ and $\lm,\mu\in \P_{m,\ell}$ with $\lm\supset \mu$.
For $u \in \Z^2$,
put
$$
\cyl\eta = \pi(\perid\lm + u),\quad \cyl\nu = \pi(\perid\mu + u).
$$
Then $\cyl\eta$ and $\cyl\nu$ are cylindric diagrams in $\Z^2/\Z(m,-\ell)$,
and
$$
f^{\cyl\lm/\cyl\mu} = f^{\cyl\eta/\cyl\nu},\ \ 
g^{\cyl\lm/\cyl\mu} = g^{\cyl\eta/\cyl\nu}.
$$
\end{prop}
%%%%%%%%%%%%%%%%%%%%%%%%%%
%

By Proposition \ref{prop:shift}, 
Theorem \ref{thm2} implies that
Conjecture \ref{conj2} is also true for $\lm=(\ell+1,1^{m-1})$ and $\mu=(0^m)$.

%%%%%%%%%%%%%%%%%%%%%%%%%%%%%%%%%%%%%%%%%%%%%%%%%%%%%
\section{Proof of Theorem \ref{thm1}} 
%%%%%%%%%%%%%%%%%%%%%%%%%%%%%%%%%%%%%%%%%%%%%%%%%%%%%%%%%%%%

Fix $\ell \geqq 1$ and $n \geqq 1$. 
Let $\lm=(n)$, $\mu=(0)$, $\cyl\lm=\cyl\lm_{(1,-\ell)}$ and
$\cyl\mu=\cyl\mu_{(1,-\ell)}$.

The cylindric skew diagram $\cy\lm/\cy\mu$ has a unique linear extension,
in which $1,2,\dots,n$ are arranged in order from  right to left.
Hence the first equality in 
Theorem \ref{thm1} holds:
$$f^{\cyl\lm/\cyl\mu}=1.$$
We will show $g^{\cyl\lm/\cyl\mu}=1$
in the rest of this section.

For $i \in \Z_{\geqq1}$,
we denote the cell $(n-i+1,1)$ by $\pt{i}$.
Note that for a subset $D$ of $\cyl\lm$,
a cell $\pt{i}$ is $D$-active if and only if
$$
\pt{i} \in D \ \hbox{and}\ 
\pt{i-1},\pt{i-\ell},\pt{i-\ell-1} \in \cy\lm\setminus D,
$$
and hence
$$
D_{\pt{i}} = \left(D \setminus \{\pt{i}\}\right) \cup \{\pt{i-\ell-1}\}.
$$
\begin{center}
\begin{tikzpicture}[scale=0.9]

\draw (3,0) -| (15,1) -- (3,1);
\foreach \x in {14,13,12,10,9,8,7,6,5,4}
\draw (\x,0) -- +(0,1);

\draw (3,-1) -| (8,1);
\foreach \x in {7,6,5,4}
\draw (\x,-1) -- +(0,1);

\foreach \x/\c in {12/3, 13/2, 14/1, 9/\ell-1, 8/\ell, 7/\ell+1, 6/\ell+2, 5/\ell+3, 4/\ell+4, 3/\ell+5}
\node at (\x+0.5,0.5) {$\pt{\c}$};
\draw[dotted, ultra thick] (2.5,0.5) -- +(-0.5,0);

\draw[dotted, ultra thick] (10.5,0.5) -- +(1,0);

\foreach \x/\c in {10/5, 11/4, 12/3, 13/2, 14/1}
\node at (\x-6.5,-0.5) {$\pt{\c}$};
\draw[dotted, ultra thick] (2.5,-0.5) -- +(-0.5,0);

\draw[dotted, ultra thick] (5.5,-1.5) -- +(0,-0.5);

\end{tikzpicture}
\end{center}

Take $q,r \in \Z$ such that
\begin{equation}\label{eq:qr}
n=q(\ell+1)+r,\quad q\geqq0,\quad 0\leqq r \leqq \ell.
\end{equation}

Put
$$
\D_{\ell;n}:= \{
(i_1,\dots,i_q) \in \Z^q \mid
i_1\geqq r+1,\ 
i_{k+1}-i_k \geqq \ell+1 \ (1 \leqq k \leqq q-1)
\}
$$
For $(i_1,\dots,i_q)\in \D_{\ell;n}$, define
$$\psi(i_1,\dots,i_q) = 
\cyl\lm \setminus \left([\pt{1},\pt{r}] \cup
\left( \bigcup_{k=1}^q [\pt{i_k},\pt{i_k+\ell}] \right) \right),
$$
where $[\pt{i},\pt{j}]=\{\pt{i},\pt{i+1},\dots,\pt{j}\}\subset \cy\lm$ for $i \leqq j$
and $[\pt{i},\pt{j}]=\emptyset$ for $i>j$.
%%%%%%%%%%%%%%%%%%%%%%%%%%%%%%%%%%%%%%%%%%%%%%%%%%%%%%%%%%%
\begin{prop}\label{prop104}
The map $\psi$
gives a bijection from $\D_{\ell;n}$ to $\E_\cyl\lm(\cyl\mu)$.
\end{prop}
%%%%%%%%%%%%%%%%%%%%%%%%%%%%%%%%%%%%%%%%%%%%%%%%%%%%%%%%%%%%%%
\begin{pf}
First,
we show that $\psi(i_1,\dots,i_q) \in \E_\cyl\lm(\cyl\mu)$.
We proceed by induction on
$$
M=\sum_{k=1}^q i_k.
$$
The number $M$ takes the minimum value 
$$M_\mathrm{min} := \frac{1}{2}q(q-1)(\ell+1)+q(r+1)$$ 
when $i_1=r+1$, $i_2=(\ell+1)+r+1$,$\dots$, $i_q=(q-1)(\ell+1)+r+1$.
For such $(i_1,\dots,i_q)$, we have
\begin{align*}
\psi(i_1,\dots,i_q) = \cyl\mu \in \E_\cyl\lm(\cyl\mu).
\end{align*}
Let $M>M_{\mathrm{min}}$
and suppose that $\psi(i_1,\dots,i_q) \in \E_\cyl\lm(\cyl\mu)$
for all $(i_1,\dots,i_q) \in \D_{\ell;n}$ such that $\sum_{k=1}^q i_k \leqq M-1$.
Take $\boldsymbol{i}=(i_1,\dots,i_q) \in \D_{\ell;n}$ with $\sum_{k=1}^q i_k =M$.
As $M>M_\mathrm{min}$,
there exists $g \in [2,q+1]$ such that $i_g-i_{g-1}>\ell+1$ (we consider ``$i_{q+1}=+\infty$'').
For such $g$,
we have $(i_1,\dots,i_{g}-1,\dots,i_q) \in \D_{\ell;n}$.
Put $D=\psi(i_1,\dots,i_g-1,\dots,i_q)$.
By induction hypothesis, $D \in \E_\cyl\lm(\cyl\mu)$.
Now $\pt{i_g+\ell} \in D$.
Note that $[\pt{i_g-1},\pt{i_g-1+\ell}] \subset \cyl\lm\setminus D$.
In particular, $\pt{i_g-1},\pt{i_g},\pt{i_g+\ell-1} \in \cyl\lm \setminus D$,
and hence the cell $\pt{i_g+\ell}$ is $D$-active.
Hence
\begin{align*}
D_{\pt{i_g+\ell}}=\left(
D\setminus \{\pt{i_g+\ell}\}\right)
\cup \{ \pt{i_g-1}\} = \psi(i_1,\dots,i_g,\dots,i_q) \in \E_\cyl\lm(\cyl\mu).
\end{align*}

Next, we show that the map $\psi$ is surjective (injectivity is obvious).
It is obvious that $\cyl\mu=\psi(i_1,\dots,i_q)$ as $i_k=(k-1)(\ell+1)+r$.
Take $D \in \E_\cyl\lm(\cyl\mu)$. 
Suppose that there exists $(i_1,\dots,i_q) \in \D_{\ell;n}$ such that
$D=\psi(i_1,\dots,i_q)$.
Any $D$-active cell is of the form $\pt{i_g+\ell+1}$
for some $g \in [1,q]$ such that $i_{g+1}-i_g>\ell+1$. 
Since $i_{g+1}-i_g>\ell+1$, $(i_1,\dots,i_g+1,\dots,i_q) \in \D_{\ell;n}$.
We have
$$
D_{\pt{i_g+\ell+1}} = \psi(i_1,\dots,i_g+1,\dots,i_q) \in \E_\cyl\lm(\cyl\mu).
$$
Hence $\psi$ is surjective.
\qed
\end{pf}

We will use the following lemma, which was suggested by H.~Tagawa.
%%%%%%%%%%%%%%%%%%%%%%%%%%%%%%%%%%%%%%%%%%%%%%%%%%%%%%%%%%%%%%%%%%%%%%%%%%%%%%%%
\begin{lem}\label{prop101}
Let $\ell,c,q\in\Z_{\geqq1}$, $r \in \Z_{\geqq0}$ and
let $(a_i)_{i\geqq1}$ be a numerical sequence such that 
\begin{itemize}
\item
$a_i \neq 0$ for all $i \geqq 1$,

\item
$\lim_{i \to \infty} a_i = +\infty$,

\item
$a_{i+\ell}-a_i=c$ for all $i\geqq 1$. 
\end{itemize}
Then
\begin{equation}\label{eq101}
\sum_{
\substack{(i_1,\dots,i_q); \\ i_1 \geqq r+1 \\ i_{k+1}-i_k \geqq \ell+1 \ (k=1,2,\dots,q-1)}}
\frac{1}{\prod_{v=1}^q \prod_{u=0}^{\ell} a_{i_v+u}}
=\frac{1}{q!c^q \prod_{u=0}^{\ell q-1} a_{r+1+u}}.
\end{equation}
\end{lem}
%%%%%%%%%%%%%%%%%%%%%%%%%%%%%%%%%%%%%%%%%%%%%%%%%%%%%%%%%%%%%%%%%%%%%%%%%%%%%%%%

\begin{pf}
By assumption, for $j,t\in\Z_{\geqq1}$, we have
\begin{equation}\label{eq102}
\frac{1}{\prod_{u=0}^{\ell t-1} a_{j+u}} - \frac{1}{\prod_{u=1}^{\ell t} a_{j+u}}
=\frac{a_{j+\ell t} - a_j}{\prod_{u=0}^{\ell t} a_{j+u}}
=\frac{ct}{\prod_{u=0}^{\ell t} a_{j+u}}.
\end{equation}
We proceed by induction on $q$ to prove (\ref{eq101}).
If $q=1$,
then we have
\begin{align*}
\sum_{i_1=r+1}^\infty \frac{1}{\prod_{u=0}^{\ell} a_{i_1+u}} 
&= \frac{1}{c} \sum_{i_1=r+1}^\infty \frac{c}{\prod_{u=0}^{\ell} a_{i_1+u}} \\
&= \frac{1}{c} \sum_{i_1=r+1}^\infty 
\left( \frac{1}{\prod_{u=0}^{\ell-1} a_{i_1+u}} - \frac{1}{{\prod_{u=1}^{\ell} a_{i_1+u}}} \right) 
\tag{by Eq. (\ref{eq102})}\\
&= \frac{1}{c} \sum_{i_1=r+1}^\infty
\left( \frac{1}{\prod_{u=0}^{\ell-1} a_{i_1+u}} - \frac{1}{{\prod_{u=0}^{\ell-1} a_{i_1+1+u}}} \right) \\
&= \frac{1}{c} \lim_{N \to \infty} \sum_{i_1=r+1}^N
\left( \frac{1}{\prod_{u=0}^{\ell-1} a_{i_1+u}} - \frac{1}{{\prod_{u=0}^{\ell-1} a_{i_1+1+u}}} \right) \\
&= \frac{1}{c \prod_{u=0}^{\ell-1} a_{r+1+u}} - \lim_{N \to \infty} \frac{1}{c \prod_{u=0}^{\ell-1} a_{N+1+u}} \\
&=  \frac{1}{c \prod_{u=0}^{\ell-1} a_{r+1+u}}.
\end{align*}

Suppose $q>1$. %and the equation \eqref{eq101} holds for smaller number.
We have
\begin{align*}
& \sum_{\substack{(i_1,\dots,i_q); \\ i_1 \geqq r+1 \\ i_{k+1}-i_k \geqq \ell+1 
\ (k=1,2,\dots,q-1)}}
\frac{1}{\prod_{v=1}^q \prod_{u=0}^{\ell} a_{i_v+u}} \\
&= \sum_{i_1=r+1}^\infty \frac{1}{\prod_{u=0}^{\ell} a_{i_1+u}}
\left( \sum_{\substack{(i_2,\dots,i_q); \\ i_2 \geqq i_1+\ell+1 \\ 
i_{k+1}-i_k \geqq \ell+1 \ (k=2,3,\dots q-1)}}
\frac{1}{\prod_{v=2}^q  \prod_{u=0}^{\ell} a_{i_v+u}} \right) \\
&= \sum_{i_1=r+1}^\infty
\frac{1}{\prod_{u=0}^{\ell} a_{i_1+u}} \cdot \frac{1}{(q-1)!c^{q-1} \prod_{u=0}^{\ell (q-1)-1} a_{i_1+\ell+1+u}} \tag{by induction hypothesis} \\
&= \frac{1}{(q-1)!c^{q-1}}
\sum_{i_1=r+1}^\infty \frac{1}{\prod_{u=0}^{\ell q} a_{i_1+u}} \\
&= \frac{1}{q!c^{q}} \sum_{i_1=r+1}^\infty
\left( \frac{1}{\prod_{u=0}^{\ell q-1} a_{i_1+u}} - \frac{1}{\prod_{u=1}^{\ell q} a_{i_1+u}} \right) 
\tag{by  Eq. (\ref{eq102})} \\
&= \frac{1}{q!c^{q} \prod_{u=0}^{\ell q-1} a_{r+1+u}},
\end{align*}
and hence we completed the proof of the (\ref{eq101}).
\qed
\end{pf}
%%%%%%%%%%%%%%%%%%%%%%%%%

The hook length for each cell $\pt{i}$ is
$$
h_\cyl\lm(\pt{\ell t+j}) = (\ell+1)t+j \ \text{ for }\  t\in\Z_{\geqq0},\ j=1,2,\dots,\ell.
$$
Put $h_i := h_\cyl\lm(\pt{i})$.
By  Proposition \ref{prop104},
we have
\begin{equation*}\label{eq5}
g^{\cyl\lm/\cyl\mu} = \frac{n!}{\prod_{i=1}^r h_i} \times
\sum_{(i_1,\dots,i_q)\in \D_{\ell;n}}
\frac{1}{\prod_{v=1}^q \prod_{u=0}^\ell h_{i_v+u}}.
\end{equation*}
By applying Lemma \ref{prop101} with
$a_i = h_i$, $c=\ell+1$ and $q,r$ as in \eqref{eq:qr},
we have
\begin{align*}
n! / g^{\cyl\lm/\cyl\mu}
&= q!(\ell+1)^q \prod_{u=1}^{q\ell+r}h_u \\
&= \left(\prod_{k=0}^{q-1} \left((k+1)(\ell+1)  \prod_{j=1}^\ell h_{k\ell+j} \right)\right) h_{q\ell+1}h_{q\ell+2}\cdots h_{q\ell+r} \\
&= \prod_{k=0}^{q-1} \Biggl( ((\ell+1)k+1)((\ell+1)k+2)\cdots((\ell+1)k+\ell-1)((\ell+1)(k+1)) \Biggr) \\
&\qquad\qquad \times (q(\ell+1)+1)(q(\ell+1)+2)\cdots(q(\ell+1)+r) \\
&= (q(\ell+1))!(q(\ell+1)+1)(q(\ell+1)+2)\cdots(q(\ell+1)+r) \\
&= (q(\ell+1)+r)! =n!.
\end{align*}
Therefore $g^{\cyl\lm/\cyl\mu}=1$ and
Theorem \ref{thm1} has been proved.
%\qed

%%%%%%%%%%%%%%%%%%%%%%%%%%%%%%%%%%%%%%%%%%%%%%%%%%%%%%%%%%%%%%%%%%%%%%%%%%%%%%%%
\section{Proof of Theorem \ref{thm2}}
%%%%%%%%%%%%%%%%%%%%%%%%%%%%%%%%%%%%%%%%%%%%%%%%%%%%%%%%%%%%%%%%%%%%%%%%%%%%%%%%

Let $\ell,m \in\Z_{\geqq1}$ and
let $\lambda=((\ell+1)^m)$, $\mu=(\ell^{m-1},0)$, $\cyl\lm=\cyl\lm_{(m,-\ell)}$
and $\cyl\mu=\cyl\mu_{(m,-\ell)}$.
Put $n=|\lm/\mu|=\ell+m$.
Any linear extension $\vep$ of $\cy\lm/\cy\mu$ satisfies
$$
\vep^{-1}(1) = (m,\ell+1), \quad
\vep^{-1}(n) = (m,1).
$$
%Remark that $n=\ell+m$.
Hence $\vep$ is uniquely determined by choosing $(m-1)$ vertical components from $[2,n-1]$.
Therefore,  the first equality in Theorem \ref{thm2} holds:
$$
f^{\cy\lm/\cy\mu}=\binom{n-2}{m-1}=\binom{\ell+m-2}{m-1}.
$$

We will prove the second equality 
$$
g^{\cy\lm/\cy\mu}=\binom{\ell+m-2}{m-1}.
$$
in the rest of this section.

For $u,v \in \Z^2$,
we write $u \to v$ if $v-u=(1,0)$ or $(0,-1)$.

%%%%%%%%%%%%%%%%%%%%%%%%%%%%%%%%%%%%%%%%%%%%%%%%%%%%%%%%%%%%%%%%%%%%%%
\begin{df}
 For $u=(a,b), v=(c,d) \in \Z^2$ with $a<c$ and $b>d$,
 a subset
$$
\path=\{u=u_1,u_2,\dots,u_r=v\}
$$
of $\Z^2$
is called a \alert{lattice path} from $u$ to $v$ 
if
$$
u_1 \to u_2\to \dots\to u_r,
$$
%for all $i=1,\dots,r-1$,
and the whole set is denoted by $\LP(u,v)$.
\if0
\smallskip\noindent
(2) Let $\theta$ be a subset os $\Z^2$.
If a lattice path $\path\in \LP(u,v)$ is a subset of $S$, then $\path$ is 
called a lattice path in $\theta$, and the whole set is denoted by
$\LP_\theta(u,v)$.
\fi
\end{df}
%%%%%%%%%%%%%%%%%%%%%%%%%%%%%%%%%%%%%%%%%%%%%%%%%%%%%%%%%%%%%%%%%%%%%%%%%%

Let $m,\ell\in\Z_{\geqq 1}$ and
let $\pi$  denote
the natural projection $\Z^2\to \Z^2/\Z(m,-\ell)=\mathcal{C}_{(m,-\ell)}$ as before.
For $u,v \in \mathcal{C}_{(m,-\ell)}$,
we write $u \to v$ if there exist
$\tilde{u}\in\pi^{-1}(u)$ and $\tilde{v}\in\pi^{-1}(v)$
such that $\tilde{u}\to \tilde{v}$.

%%%%%%%%%%%%%%%%%%%%%%%%%%%%%%%%%%%%%%%%%%%%%%%%%%%%%%%%%%%%%%%%%%%%%
\begin{df}
\smallskip\noindent
 A subset
$$
\path=\{u_1,u_2,\dots,u_n\}
$$
of $\mathcal{C}_{(m,-\ell)}$
is called a \alert{non-intersecting loop} in $\mathcal{C}_{(m,-\ell)}$ 
if $n=\ell+m$ and
$$
u_1 \to u_{2}\to\dots\to u_{n-1}\to u_n\to u_1.
$$
%and 
%$$i\neq j\Rightarrow u_i\neq u_j$$. 
The whole set of non-intersecting loops is denoted by $\LL$.
\end{df}
%%%%%%%%%%%%%%%%%%%%%%%%%%%%%%%%%%%%%%%%%%%%%%%%%%%%%%%%%%%%%%%%%%%%%%%%%%

%
%Let $\lambda=((\ell+1)^m)$, $\mu=(\ell^{m-1},0) \in \P_{m,\ell}$,
%and

Let $\lm=(\lm_1,\dots,\lm_m)
\in \P_{m,\ell}$ and let $\semid{\lm}$ denote the semi-infinite diagram
$$\semid{\lm} = \{ (a,b) \in \mathbb{Z}^2 \mid 1 \leqq a \leqq m,\ b \leqq \lm_a \}$$
as before.
Note that $\semid\lambda$
is in one-to-one correspondence
with the cylindric diagram $\cyl{\lm}$ 
via the projection $\pi$.

%For $u,v\in\semid{\lm}$, 
Define
\begin{align*}
\LP_{{\lm}}(u,v)&=\{\path\in \LP(u,v)\mid \path\subset \semid{\lm}\}
\ \ (u,v\in \semid{\lm}),
\\
\LL_{{\lm}}&=\{\path\in \LL\mid \path\subset \cyl{\lm}\}.
\end{align*}

%In this case, $\LP_{{\lm}}(u,v)=\LP(u,v)$ for $u,v\in\semid{\lm}$.
% $\cyl{\lm}$ via the projection $\pi:\Z^2\to\mathcal{C}_{(m,-\ell)}$.
%%%%%%%%%%%%%%%%%%%%%%%%%%%%%%%%%%%%%%%%%%%%%%%%%%%%%%%%%%%%%%
\begin{lem}\label{lem;PtoP}
Let $\lm=(\lm_1,\dots,\lm_m)\in \P_{m,\ell}$.
Then the projection $\pi$ induces a bijection 
$$
\bigsqcup_{i=0}^\infty \LP_{{\lm}}((1,\lm_1-i),(m,\lm_1-\ell-i))
\stackrel{\cong}{\longrightarrow} \LL_\lm.
$$
\end{lem}
%%%%%%%%%%%%%%%%%%%%%%%%%%%%%%%%%%%%%%%%%%%%%%%%%%%%%%%%%%%%%%%
\begin{pf}
For  a lattice path $\path\in \LP_{{\lm}}((1,\lm_1-i),(m,\lm_1-\ell-i))$,
it is clear that  $\pi(\path)\in\LL_\lm$.
The inverse map is given by $\path\mapsto
\pi^{-1}(\path)\cap \semid{\lm}$.
\qed
\end{pf}
%%%%%%%%%%%%%%%%%%%%%%%%%%%%%%%%%%%%%%%%%%%%%%%%%%%%%%%%%%%%%%%%%%%%%
Now, we return to the special case where $\lm=((\ell+1)^m)$.
In this case, we have
\begin{equation}\label{eq;LP=LP}
\LP_{{\lm}}(u,v)=\LP(u,v)
\end{equation}
for all $u,v\in\semid{\lm}$.
%%%%%%%%%%%%%%%%%%%%%%%%%%%%%%%%%%%%%%%%%%%%%%%%%%%%%%%%%%
\begin{prop}\label{prop;EtoP}
Let $\ell,m\in\Z_{\geqq1}$ and $\lm=((\ell+1)^m)\in  \P_{m,\ell}$.

\smallskip\noindent
$(1)$ Let $k\in[0,\ell]$ and let $\nu^{(k)}=(\ell^{m-1},k)\in \P_{m,\ell}$.
Then the correspondence $\path \mapsto {\lm} \setminus \path$ gives a bijection 
$$\LP((1,\ell+1),(m,k+1))\stackrel{\cong}{\longrightarrow}
\E_{{\lm}}({\nu^{(k)}}) .$$
%
%\smallskip
\noindent
$(2)$ Let $\mu=(\ell^{m-1},0)\in \P_{m,\ell}$.
Then the correspondence $\path \mapsto \cyl{\lm} \setminus \path$ gives a bijection 
$$\LL_\lm \stackrel{\cong}{\longrightarrow} \E_{\cyl{\lm}}(\cyl{\mu}).$$
%where $\bar{D} = \pi^{-1}(D) \cap \semid{\lm}.$
%([1,m] \times \Z)$.
\end{prop}
%%%%%%%%%%%%%%%%%%%%%%%%%%%%%%%%%%%%%%%%%%%%%%%%%%%%%%%%%%%%%%%%%%%%%%
\begin{pf}
(1) 
We write $\LP=\LP((1,\ell+1),(m,k+1))$ and $\nu=\nu^{(k)}$ in this proof.

First, we will prove $\psi(\path): ={\lm}\setminus\path$ is contained in
$\E_{{\lm}}(\nu)$ by 
induction on the number $r(\path)$ of cells in $\lm$ 
which is located to the right of $\path$.
If $r(\path)=0$ then $\psi(\path)=\nu\in  \E_{{\lm}}({\nu})$.

Let $\path\in\LP$ with $r(\path)\geqq 1$
and suppose that $\psi(\mathfrak{q})\in \E_{{\lm}}({\nu})$
for any $\mathfrak{q}\in \LP$ such that $r(\mathfrak{q})< r(\path)$.
There exists
$y=(a,b)\in \lm\setminus\path$ such that
$$(a,b-1),(a-1,b), x=(a-1,b-1)\in \path.$$

Now $\path':=(\path\setminus\{x\})\cup\{y\}$ is a lattice path with
$r(\path')=r(\path)-1$ and hence $\psi(\path')\in\E_{{\lm}}({\nu})$.
Moreover, $\psi(\path)$ is obtained from $\psi(\path')$ by 
applying the elementary excitation at $x$.
This implies $\psi(\path)\in\E_{{\lm}}({\nu})$.

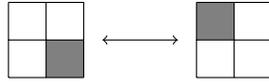
\begin{center}
\begin{tikzpicture}[scale=0.5]
\fill[gray] (1,0) rectangle +(1,1);
\draw (0,0) grid (2,2);
%\node at (0.5,1.5) {$y$};
%\node at (1.5,0.5) {$x$};

\draw[<->] (2.5,1) -- +(2,0);
%\node[above] at (7.5,1) {elementary excitation at $y$};

\begin{scope}[xshift=5cm]
\fill[gray] (0,1) rectangle +(1,1);
\draw (0,0) grid (2,2);
%\node at (0.5,1.5) {$y$};
%\node at (1.5,0.5) {$x$};
\end{scope}

\end{tikzpicture}
\end{center}

Next, we construct an inverse map.
We define $\varphi(D)=\semid{\lm} \setminus {D}$ 
and  will prove that $\varphi(D)$
is contained in $\LP$ for all $D\in \E_{{\lm}}({\nu})$.
% for all $D\in \E_{\cyl{\lm}}(\cyl{\mu})$.
For $D=\nu$, we have $\varphi(D)={\lm}/{\nu}\in \LP$.
%we have
%$$
%\varphi(\cyl\mu)=\semid\lm/\semid\mu = \{(1,\ell+1),(2,\ell+1),\dots,(m,\ell+1),(m,\ell),\dots,(m,1)\} \in \LP_0.
%$$
%
It is easy to see that 
$\varphi(D_y)\in \LP$
for any $\varphi(D)$ and any $D$-active cell $y$. 
Hence $\varphi(D)\in \LP$ for any $D\in \E_{{\lm}}({\nu})$.
This completes the proof of (1).

The statement (2) is proved by a parallel argument.
\qed
\end{pf}

Combinning \eqref{eq;LP=LP}, Lemma \ref{lem;PtoP} (2) and Proposition \ref{prop;EtoP},
we have the following:
%%%%%%%%%%%%%%%%%%%%%%%%%%%%%%%%%%%%%%%%%%%%%%%%%%%%%%%%%%%%%%%%
\begin{cor}\label{cor;EtoL}
Let $\ell,m\in \Z_{\geqq1}$ and $\lm=((\ell+1)^m),
\mu=(\ell^{m-1},0)\in \P_{m,\ell}$.
Then the map $\path\mapsto \cyl{\lm}\setminus \pi(\path)$
gives a bijection
$$
\bigsqcup_{i=0}^\infty \LP((1,\ell+1-i),(m,1-i))
\stackrel{\cong}{\longrightarrow} 
\E_{\cyl{\lm}}(\cyl{\mu}). $$
\end{cor}
%%%%%%%%%%%%%%%%%%%%%%%%%%%%%%%%%%%%%%%%%%%%%%%%%%%%%%%%%%%%%

For $\ell,m\in\Z_{\geqq1}$ and $x=(a,b) \in \Z^2$,
define 
%$h_\lambda(a,b)$ is given by
\begin{equation}\label{eq401}
h_{m,\ell}(x) = \ell+m-a-b+2.
\end{equation}
%$\lambda$のcell ???
%
Remark that 
$$h_{m,\ell}(x) =h_{\lm}(x)$$ for $\lm=((\ell+1)^m)$ and $x\in {\lm}$.
For $s \in \Z_{\geqq1}$, define
\begin{equation}\label{eq402}
F_{(\ell,m;s)} = \sum_{\path \in \LP((1,\ell+1),(m,2))} 
\prod_{x \in \path} \frac{1}{h_{m,\ell}(x)+s-1}.
\end{equation}

%%%%%%%%%%%%%%%%%%%%%%%%%%%%%%%%%%%%%%%%%%%%%%%%%%%%%%%%%%%%%%%%%%%%%%%%%%
\begin{lem}\label{lem401}
Let $\ell,m,s\in\Z_{\geqq1}$. Then 
\begin{equation}\label{eq;recFs}
F_{(\ell,m;s)} = \frac{(s-1)!}{(\ell+m+s-2)!}\binom{\ell+m-2}{m-1}.
\end{equation}
\end{lem}
%%%%%%%%%%%%%%%%%%%%%%%%%%%%%%%%%%%%%%%%%
\begin{pf}
We proceed by induction on $s$.
If $s=1$,
then it follows from Proposition \ref{prop;EtoP} (1) with $k=1$ and Theorem \ref{thm;NHF}
that
$$
F_{(\ell,m;1)}
= \sum_{\path \in \LP((1,\ell+1),(m,2))} \prod_{x \in \path} \frac{1}{h_{m,\ell}(x)}
= \frac{1}{(\ell+m-1)!}\binom{\ell+m-2}{m-1},
$$
 %
%for $\lambda/\nu$,
and hence in this case \eqref{eq;recFs} holds for all $\ell,m \in \Z_{\geqq1}$.

Take $s\geqq 1$ and 
suppose that \eqref{eq;recFs} holds  for all $\ell,m \in \Z_{\geqq1}$.
%Suppose $s>1$.
Via the bijection
$$
\LP((1,\ell+2),(m,2)) \cong 
\LP((2,\ell+2),(m,2)) \sqcup \LP((1,\ell+1),(m,2)),
$$
%where $u'$ and $u''$ satisfy $u'-u=(1,0)$ and $u''-u=(0,-1)$, respectively. Therefore, 
we have
$$
F_{(\ell+1,m;s)} = \frac{1}{h_{\ell+1,m}(1,\ell+2)+s-1} 
(F_{(\ell+1,m-1;s)} + F_{(\ell,m;s+1)}).
$$
Hence we have
\begin{align*}
F_{(\ell,m;s+1)}
&= (h_{\ell+1,m}(1,\ell+2)+s-1)
 F_{(\ell+1,m;s)} - F_{(\ell+1,m-1;s)} \\
&= (m+s-1) \cdot \frac{(s-1)!\binom{\ell+m-1}{m-1}}{(\ell+m+s-1)!}
- \frac{(s-1)!\binom{\ell+m-2}{m-2}}{(\ell+m+s-2)!} 
\tag{by induction hypothesis} \\
&= \frac{(s-1)!\binom{\ell+m-2}{m-2}}{(\ell+m+s-2)!}
\left(
\frac{(m+s-1)(\ell+m-1)}{(\ell+m+s-1)(m-1)} - 1
\right) \\
&= \frac{(s-1)!\binom{\ell+m-2}{m-2}}{(\ell+m+s-2)!}
\cdot \frac{(m+s-1)(\ell+m-1) - (\ell+m+s-1)(m-1)}{(\ell+m+s-1)(m-1)} \\
&= \frac{(s-1)!\binom{\ell+m-2}{m-2}}{(\ell+m+s-2)!}
\cdot \frac{\ell s}{(\ell+m+s-1)(m-1)} \\
&= \frac{s!}{(\ell+m+s-1)!}\binom{\ell+m-2}{m-1}.
\end{align*}
This completes the induction step. 
\qed
\end{pf}

%%%%%%%%%%%%%%%%%%%%%%%%%%%%%%%%%%%%%%%%
\if0
For $\ell,m,s,t \in \Z_{\geqq1}$ and $x\in\Z^2$,
%For $a,b,c \in \Z$ with $1 \leqq a \leqq m$ and $2 \leqq b \leqq \ell+1$ and $c \geqq 0$,
we define
\begin{equation}\label{eq;cylhook}
h^{s,t}_{m,\ell}(x)
%(a,b-c\ell) :
= \ell+m-a-d+ct+s+1,
\end{equation}
where $x=(a,b)$ and $b=c(\ell+1)+d$ with $c\geqq 0,\ 1\leqq d\leqq \ell+1$.
\fi
%%%%%%%%%%%%%%%%%%%%%%%%%%%%%%%%%%%%%%%%
For $\ell,m,s,t \in \Z_{\geqq1}$ and $x\in\Z^2$,
%For $a,b,c \in \Z$ with $1 \leqq a \leqq m$ and $2 \leqq b \leqq \ell+1$ and $c \geqq 0$,
we define
\begin{equation}\label{eq;cylhook}
h^{s,t}_{m,\ell}(x)
%(a,b-c\ell) :
= \ell+m-a-b+(d-c)t+s+1,
\end{equation}
where $x=(a+cm,b-d\ell)$ with
$1 \leqq a \leqq m,\ 2 \leqq b \leqq \ell+1$ and $c,d\geqq1$.
Remark that for $s=1$ and $t=\ell+m$, 
the number $h^{s,t}_{m,\ell}(x)$ gives a cylindric hook length:
\begin{equation}\label{eq7}
h^{1,\ell+m}_{m,\ell}(x) =  h_\cyl{\lm}(\pi(x))
\end{equation}
 for $\lm=((\ell+1)^m)$ and $x\in \semid{\lm}$.
(See Definition \ref{df;chook}.)
% for $x \in \cyl{\lm}$.

Note also that 
for any $s,t,a \in\Z_{\geqq1}$, we have
\begin{equation}\label{eq305}
\lim_{b \to -\infty} h^{s,t}_{m,\ell}(a,b) = +\infty.
\end{equation}

%For $s\in \Z_{\geqq1}$, we put
%$$
%h_{m,\ell,s}^t (x) = h_{m,\ell}^t(x)+s-1
%$$
%and define
Define
$$
F_{(\ell,m;s,t)}
= \sum_{i=0}^\infty \sum_{\path\in\LP((1,\ell+1-i),(m,1-i))} \prod_{x\in\path} 
\frac{1}{h_{m,\ell}^{s,t}(x)}.
%{h^t_{m,\ell}(x)+s-1}.
$$

%%%%%%%%%%%%%%%%%%%%%%%%%%%%%%%%%%%%%%%%%%%%
\begin{lem}\label{lem;Flmstrec}
Let $\ell,m,s,t \in \Z_{\geqq1}$. Then
\begin{equation}
F_{(\ell,m;s,t)} = \frac{1}{t-m+1}
\left(F_{(\ell,m;s)} + F_{(\ell,m-1;s+1,t)} - F_{(\ell,m-1;s,t)} \right).
\end{equation}
\end{lem}
%%%%%%%%%%%%%%%%%%%%%%%%%%%%%%%%%%%%%%%%%%%%%%
\begin{pf}
For $i\geqq 0$,
put $$d_i = \ell+1-i.$$
Fix $\ell,m,s,t$ and write 
$$
h (x) = h_{m,\ell}^{s,t}(x)
$$
for  a while.
We have
$$
h(m,1-i) - h(1,\ell+1-i) = t-m+1,
$$
and hence
\begin{align*}
&F_{(\ell,m;s,t)} \\
&= \sum_{i=0}^\infty \sum_{0\leqq k_1 \leqq k_2 \leqq \cdots \leqq k_{m-1} \leqq \ell}
\frac{1}{\prod_{r=0}^{k_1} h(1,d_i-r) 
\prod_{r=k_1}^{k_2} h(2,d_i-r) \cdots \prod_{r=k_{m-1}}^{\ell} h(m,d_i-r)} \\
&= \frac{1}{t-m+1} \sum_{i=0}^\infty \\
&\quad \times
 \sum_{0\leqq k_1 \leqq k_2 \leqq \cdots \leqq k_{m-1} \leqq \ell}
\left(
\frac{1}{\prod_{r=0}^{k_1} h(1,d_i-r) \prod_{r=k_1}^{k_2} h(2,d_i-r) \cdots \prod_{r=k_{m-1}}^{\ell-1} h(m,d_i-r)} \right. \\
&\qquad\qquad\left.
- \frac{1}{\prod_{r=1}^{k_1} h(1,d_i-r) \prod_{r=k_1}^{k_2} h(2,d_i-r) \cdots \prod_{r=k_{m-1}}^{\ell} h(m,d_i-r)}
\right) \\
&= \frac{1}{t-m+1} \sum_{i=0}^\infty (A_i+B_i-C_i-D_i).
\end{align*}

\begin{figure}[h]
\begin{center}
\begin{tikzpicture}[rotate=180]

\fill[gray] (-1,0.5) rectangle +(13,3);
\fill[white] (0.5,0.5) rectangle +(4,1);
\fill[white] (3.5,1.5) rectangle +(4,1);
\fill[white] (6.5,2.5) rectangle +(4,1);

\draw (12,0.5) -- (-1,0.5) |- (12,3.5);
\draw (12,1.5) -- (-1,1.5);
\draw (12,2.5) -- (-1,2.5);

\draw (0.5,0.5) rectangle +(4,1);
\draw (3.5,1.5) rectangle +(4,1);
\draw (6.5,2.5) rectangle +(4,1);
\draw (1.5,0.5) -- +(0,1);
\draw (3.5,0.5) -- +(0,1);
\draw (4.5,1.5) -- +(0,1);
\draw (6.5,1.5) -- +(0,1);
\draw (7.5,2.5) -- +(0,1);
\draw (9.5,2.5) -- +(0,1);

\node at (1,1) {\tiny\rotatebox{45}{$(1,d_i)$}};
\draw[ultra thick, dotted] (2,1) -- +(1,0);
\node at (4,1) {\tiny\rotatebox{45}{$(1,d_i-k_1)$}};
\node at (4,2) {\tiny\rotatebox{45}{$(2,d_i-k_1)$}};
\draw[ultra thick, dotted] (5,2) -- +(1,0);
\node at (7,2) {\tiny\rotatebox{45}{$(2,d_i-k_2)$}};
\node at (7,3) {\tiny\rotatebox{45}{$(3,d_i-k_2)$}};
\draw[ultra thick, dotted] (8,3) -- +(1,0);
\node at (10,3) {\tiny\rotatebox{45}{$(3,d_i-\ell)$}};

\end{tikzpicture}
\end{center}
\caption{An excited diagram for $\cyl{\lambda}/\cyl{\mu}$ with  $\lambda=(\ell^3)$, $\mu=((\ell-1)^2,0)$.}
\end{figure}
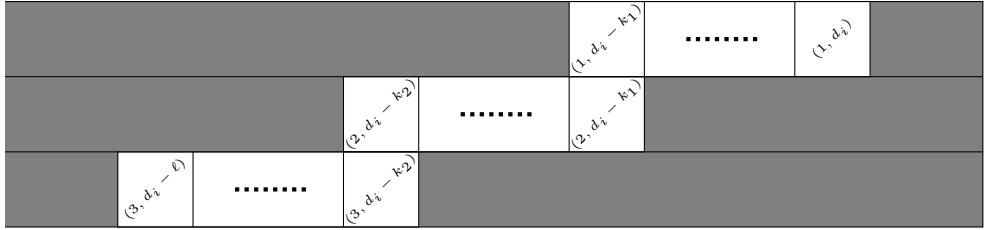

Here, $A_i$, $B_i$, $C_i$ and $D_i$ are
%in the summation are put as
\begin{align*}
A_i &= \sum_{0\leqq k_1 \leqq \cdots \leqq k_{m-1} \leqq \ell-1}
\frac{1}{\prod_{r=0}^{k_1} h(1,d_i-r) \prod_{r=k_1}^{k_2} h(2,d_i-r) 
\cdots \prod_{r=k_{m-1}}^{\ell-1} h(m,d_i-r)}, \\
B_i &= \sum_{0 \leqq k_1 \leqq k_2 \leqq \cdots \leqq k_{m-2} \leqq \ell}
\frac{1}{\prod_{r=0}^{k_1} h(1,d_i-r) \prod_{r=k_1}^{k_2} h(2,d_i-r) 
\cdots \prod_{r=k_{m-2}}^{\ell} h(m-1,d_i-r)}, \\
C_i &= \sum_{0 \leqq k_2 \leqq k_3 \leqq \cdots \leqq k_{m-1} \leqq \ell}
\frac{1}{\prod_{r=0}^{k_2} h(2,d_i-r) \prod_{r=k_1}^{k_2} h(3,d_i-r) 
\cdots \prod_{r=k_{m-1}}^{\ell} h(m,d_i-r)}, \\
D_i &= \sum_{1\leqq k_1 \leqq k_2 \leqq \cdots \leqq k_{m-1} \leqq \ell}
\frac{1}{\prod_{r=1}^{k_1} h(1,d_i-r) \prod_{r=k_1}^{k_2} h(2,d_i-r) 
\cdots \prod_{r=k_{m-1}}^{\ell} h(m,d_i-r)}\\
&=A_{i+1}.
\end{align*}

Now we have
%can write $B_i$ and $C_i$ as
\begin{align*}
\sum_{i=0}^\infty B_i
&= \sum_{i=0}^\infty \sum_{\path\in\LP((1,\ell+1-i),(m-1,1-i))}
\prod_{x\in\path} \frac{1}{h(x)} \\
&= \sum_{i=0}^\infty \sum_{\path\in\LP((1,\ell+1-i),(m-1,1-i))}
\prod_{x\in\path} \frac{1}{h^{s+1,t}_{m-1,\ell}(x)} 
%
%&=\sum_{i=0}^\infty \sum_{\path\in\LP((2,\ell+1-i),(m,1-i))}
%\prod_{x\in\path} \frac{1}{h_{s+1}^t(x)} 
%
 = F_{(\ell,m-1;s+1,t)}, \\
\sum_{i=0}^\infty C_i
&=\sum_{i=0}^\infty \sum_{\path\in\LP((2,\ell+1-i),(m,1-i))}
\prod_{x\in\path} \frac{1}{h(x)}
 = F_{(\ell,m-1;s,t)}.
\end{align*}
Moreover, using \eqref{eq305}, we have
\begin{align*}
\sum_{i=0}^\infty (A_i-D_i)
&= \sum_{i=0}^\infty (A_i - A_{i+1}) = A_0 \\
&= \sum_{i=0}^\infty \sum_{\path\in\LP((1,\ell+1),(m,2))}
\prod_{x\in\path} \frac{1}{h(x)}
 = F_{(\ell,m;s)}. 
\end{align*}
Therefore, 
\begin{align*}
&F_{(\ell,m;s,t)} \
= \frac{1}{t-m+1}
\left(F_{(\ell,m;s)} + F_{(\ell,m-1;s+1,t)} - F_{(\ell,m-1;s,t)} \right). 
\end{align*}
\qed\end{pf}

%%%%%%%%%%%%%%%%%%%%%%%%%%%%%%%%%%%%%%%%%%%%%%%%%%%%%%%%%%%%%%%%%%%%%
\begin{prop}\label{prop103}
Let $\ell,m,s,t \in \Z_{\geqq1}$. Then
\begin{equation}\label{eq;Flmst}
F_{(\ell,m:s,t)} = \frac{(s-1)!}{(\ell+m+s-2)! t}\binom{\ell+m-2}{m-1}.
\end{equation}
\end{prop}
%%%%%%%%%%%%%%%%%%%%%%%%%%%%%%%%%%%%%%%%%%%%%%%%%%%%%%%%%%%%%%%%%%%

\begin{pf}
\if0
For simplicity, we write
$$
h^t (x) = h_{\cyl\lm}^t(x)+s-1
$$
for $x \in \cyl\lm$ in this proof.
%
\fi
We proceed by induction on $m$.
%
%By \eqref{eq;cylhook},
%Since $t \geqq 1$, 
If $m=1$,
then putting $d_i=\ell+1-i$, we have
\begin{align*}
F_{(\ell,1;s,t)}
&= \sum_{i=0} \frac{1}{\prod_{k=0}^\ell h^{s,t}_{1,\ell}(1,d_i-k)} \\
&= \frac{1}{t} \sum_{i=0} 
\left(\frac{1}{\prod_{k=0}^{\ell-1} h^{s,t}_{1,\ell}(1,d_i-k)} - 
\frac{1}{\prod_{k=1}^{\ell} h^{s,t}_{1,\ell}(1,d_i-k)}\right) \\
%&= \frac{1}{t} \cdot \frac{1}{h^t_{1,\ell,s}(1,\ell+1) 
%h^t_{1,\ell,s}(1,\ell) \cdots h^t_{1,\ell,s}(1,2)} 
&= \frac{1}{t} 
\frac{1}{ \prod_{k=0}^{\ell-1} h^{s,t}_{1,\ell}(1,\ell+1-k) } 
\tag{by (\ref{eq305})} \\
&= \frac{1}{t} \cdot \frac{1}{s (s+1) \cdots (s+\ell-1)} 
\end{align*}
for any $\ell,s,t \in \Z_{\geqq1}$.
This proves
\eqref{eq;Flmst} when $m=1$.

Let $m>1$ and suppose that
\begin{equation}%\label{eq;Flmst}
F_{(\ell,m-1:s,t)} = \frac{(s-1)!}{(\ell+(m-1)+s-2)! t}\binom{\ell+(m-1)-2}{m-2}.
\end{equation}
for all $\ell,s,t$.
By Lemma \ref{lem;Flmstrec}, we have
%Let $m>1$.
\begin{align*}
F_{(\ell,m;s,t)} &= \frac{1}{t-m+1}
\left(F_{(\ell,m;s)} + F_{(\ell,m-1;s+1,t)} - F_{(\ell,m-1;s,t)} \right)\\
&= \frac{1}{t-m+1}
\left(
\frac{(s-1)!\binom{\ell+m-2}{m-1}}{(\ell+m+s-2)!}
+ \frac{s!\binom{\ell+m-3}{m-2}}{(\ell+m+s-2)!t}
- \frac{(s-1)!\binom{\ell+m-3}{m-2}}{(\ell+m+s-2)!t}
\right) 
\tag{by Lemma \ref{lem401} and induction hypothesis}\\
&= \frac{(s-1)!\binom{\ell+m-3}{m-2}}{(\ell+m+s-3)!(t-m+1)}
\left(
\frac{\ell+m-2}{(m-1)(\ell+m+s-2)}
+ \frac{s}{(\ell+m+s-2)t}
- \frac{1}{t}
\right) \\
&= \frac{(s-1)!\binom{\ell+m-3}{m-2}}{(\ell+m+s-3)!(t-m+1)}
\cdot \frac{(\ell+m-2)t + s(m-1) - (m-1)(\ell+m+s-2)}{(m-1)(\ell+m+s-2)t} \\
&= \frac{(s-1)!}{(\ell+m+s-3)!(t-m+1)}\binom{\ell+m-3}{m-2}
\cdot \frac{(\ell+m-2)(t-m+1)}{(m-1)(\ell+m+s-2)t} \\
&= \frac{(s-1)!}{(\ell+m+s-2)!t}\binom{\ell+m-2}{m-1}.
\end{align*}
We have proved Proposition \ref{prop103}.
\qed
\end{pf}

Finally,
%we compute $g^{\cyl\lm/\cyl\mu}$
%the right hand side of (\ref{eq;conj2}) 
by applying Corollary \ref{cor;EtoL} and
Proposition \ref{prop103} with $s=1$ and $t=\ell+m=n$,
%and then 
we obtain
\begin{align*}
g^{\cyl\lm/\cyl\mu}
&= n! \sum_{D\in\E_\cyl{\lm}(\cyl{\mu})} \prod_{x\in \cyl{\lm}\setminus D} 
\frac{1}{h_{\cyl\lm}(x)} \\
&= n! \sum_{i=0}^\infty \sum_{\path\in\LP((1,n-m+1-i),(m,1-i))} \prod_{x\in\path} 
\frac{1}{h^{1,\ell+m}_{m,\ell}(x)} \\
&= n! \cdot F_{(\ell,m;1,\ell+m)} \\
&= n! \cdot \frac{1}{(\ell+m-1)!(\ell+m)}\binom{\ell+m-2}{m-1} \\
&= \binom{\ell+m-2}{m-1}.
\end{align*}
This completes the proof of Theorem \ref{thm2}.
%\qed
%\bibliographystyle{abbrv}
%\bibliography{CoHFoCSD}
%\nocite{MPP,Str,Proc1,Proc2,Proc3}
%\if0

\end{document}